\documentclass[5pt,a4paper]{article}%
\usepackage{bm}
\usepackage{eurosym}
\usepackage{amsmath}
\usepackage{amsfonts}
\usepackage{amssymb}
\usepackage{graphicx}
\usepackage{fullpage}
\usepackage{url}
\usepackage{amsthm}
\usepackage{caption}
\usepackage{subcaption}
\usepackage{comment}%
\usepackage[title]{appendix}
\usepackage{cite}
\usepackage{hyperref}
\usepackage{color}

\providecommand{\U}[1]{\protect\rule{.1in}{.1in}}

\newtheorem{lemma}{Lemma}
\newtheorem{proposition}{Proposition}
\providecommand{\keywords}[1]
{
  \small	
  \textbf{\textit{Keywords---}} #1
}

\newcommand{\inti}[2][y]{\int_{0}^{\infty} #2 \, d#1}
\newcommand{\sech}{\text{sech}}
\newcommand{\dd}{\text{d}}

\numberwithin{equation}{section}

\begin{document}

\title{Stability and slow dynamics of an interior spiky pattern in a one-dimensional spatial Solow  model with capital-induced labor migration}
\author{Fanze Kong\footnote{Department of Applied Mathematics, University of Washington, Seattle, 98195, WA, USA, email address: {fzkong@uw.edu}.  }, Jiayi Sun \footnote{School of Mathematics, Hunan University, Changsha, 410082, P. R. China, email address: {sjy1205@hnu.edu.cn}} and Shuangquan Xie\footnote{School of Mathematics, Hunan University, Changsha, 410082, P. R. China, email address: {xieshuangquan@hnu.edu.cn}}.} 
\maketitle

\begin{abstract}
One of the most significant findings in the study of spatial Solow-Swan models is the emergence of economic agglomeration, in which economic activities concentrate in specific regions.  Such agglomeration provides a fundamental mechanism driving the spatial patterns of urbanization, labor migration, productivity growth, and resource allocation.  In this paper, we consider the one-dimensional spatial Solow-Swan model with capital-induced labor migration, which captures the dynamic interaction between labor and capital through migration and accumulation.  Focusing on the regime of sufficiently small capital diffusivity, we first construct an interior spike (spiky economic agglomeration) quasi-equilibrium.  Next, we perform the linear stability of the corresponding spike equilibrium by using a hybrid asymptotic and numerical method.  We show that a single interior spike remains stable for small reaction-time constants but undergoes a Hopf bifurcation when the constant is sufficiently large, leading to oscillations in spike height (economic fluctuation). Finally, we derive a differential–algebraic system to capture the slow drift motion of quasi-equilibrium (core-periphery shift).  Numerical simulations are carried out to support our theoretical studies and reveal some intriguing yet unexplained dynamics.  
\end{abstract}
\section{Introduction}
This paper is devoted to the following reaction-diffusion-advection system of the labor-capital pair $(l,k)$, which models spatial economic activities by capturing how capital accumulation drives labor migration subject to the labor growth: 
\begin{align}\label{tt1}
\left\{
\begin{array}{ll}
l_t =\overbrace{d_l l_{xx}}^{\text{labor diffusion}} - \overbrace{\chi (l k_x)_x}^{\text{captial-induced flux}} + \overbrace{l(a - bl)}^{\text{labor population source}}, & x \in (-\ell_b,\ell_b),~ t > 0, \\
k_t =\overbrace{ d_k k_{xx}}^{\text{capital diffusion}} - \overbrace{\delta k}^{\text{capital depreciation}}+\overbrace{ (1-c)k^\theta l^{1-\theta}}^{\text{Cobb–Douglas production function}}, & x \in (-\ell_b,\ell_b),~ t > 0, \\
l_x(x,t) = k_x(x,t) = 0, & x = 0,l, ~t > 0, \\
l(x,0) = l_0(x) \geq \not\equiv 0, k(x,0) = k_0(x) \geq \not\equiv 0, & x \in (-\ell_b,\ell_b),
\end{array}
\right.
\end{align}
where $x$ and $t$ are spatial and temporal variables; $\ell _b>0$ is a constant; $l(x,t)$ and $k(x,t)$ denote the labor population density 
and capital distribution, respectively. Here positive constants $d_l$, $d_k$ represent the labor and capital diffusion rates, respectively; costant $\chi>0$ measures the strength of labor migration in response to variations in capital concentration; $a$ denotes the intrinsic growth rate of labor; $b$ characterizes the strength of intraspecific competition. The parameter $\theta \in (0,1)$ represents capital elasticity of output
, while $c\in(0,1)$ denotes the savings/consumption rate, leaving $(1-c)$ for reinvestment.  In particular, the Neumann boundary conditions ($l_x = k_x = 0$) reflect our assumption of a closed economy in region $\Omega$, where neither labor nor capital crosses the boundary.

\begin{figure}
    \centering
    \begin{subfigure}[b]{0.34\linewidth}
        \centering
        \includegraphics[width=\linewidth]{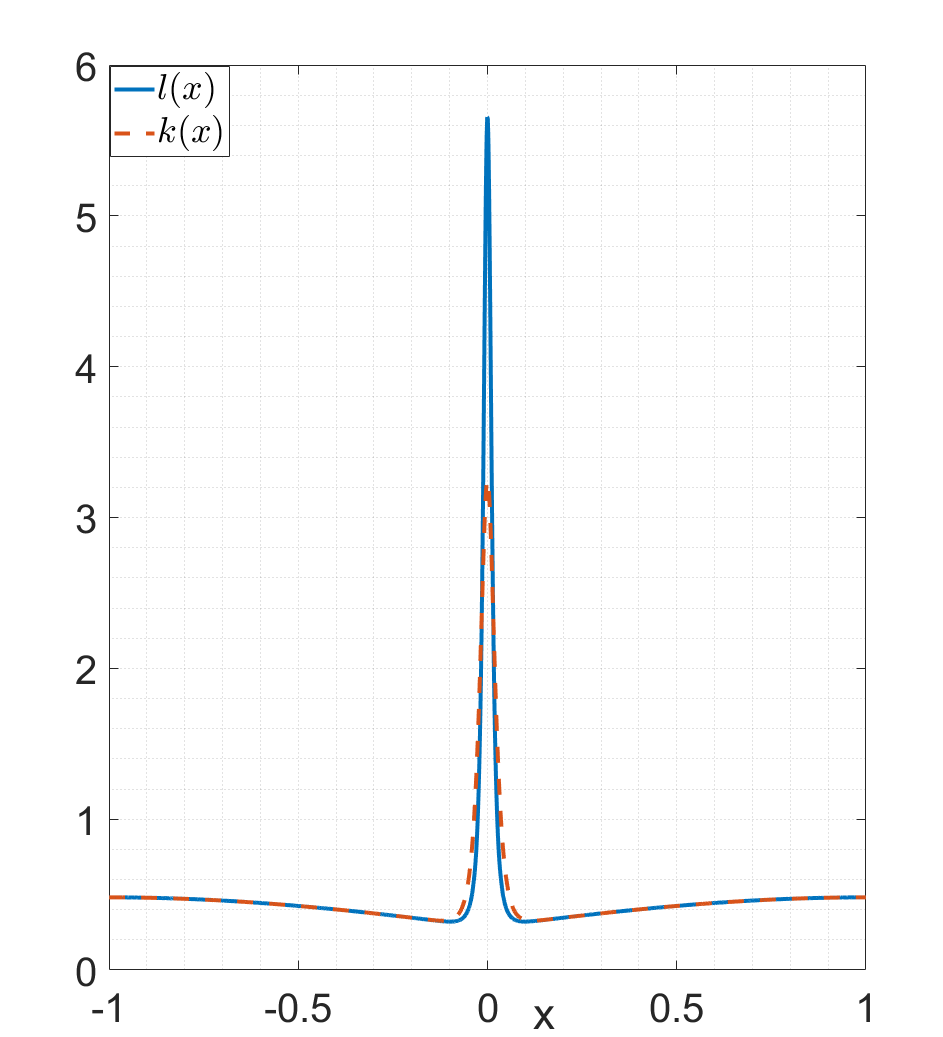}
        \caption{Single spike profile at steady state.}\label{fig:intro(a)}
    \end{subfigure}
    \hfill
    \begin{subfigure}[b]{0.32\linewidth}
        \centering
        \includegraphics[width=\linewidth]{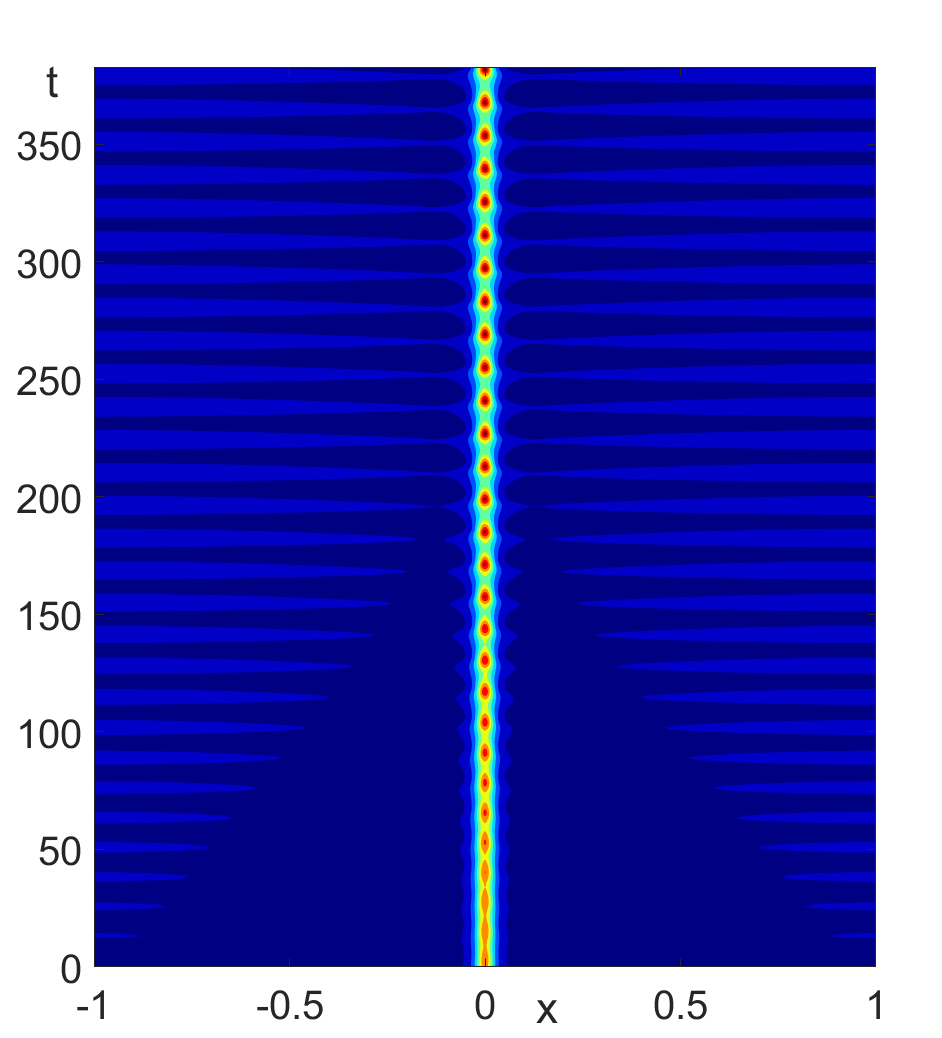}
        \caption{Height oscillation of a single spike ($\tau=2.7$).}\label{fig:intro(b)}
    \end{subfigure}
    \hfill
    \begin{subfigure}[b]{0.32\linewidth}
        \centering
        \includegraphics[width=\linewidth]{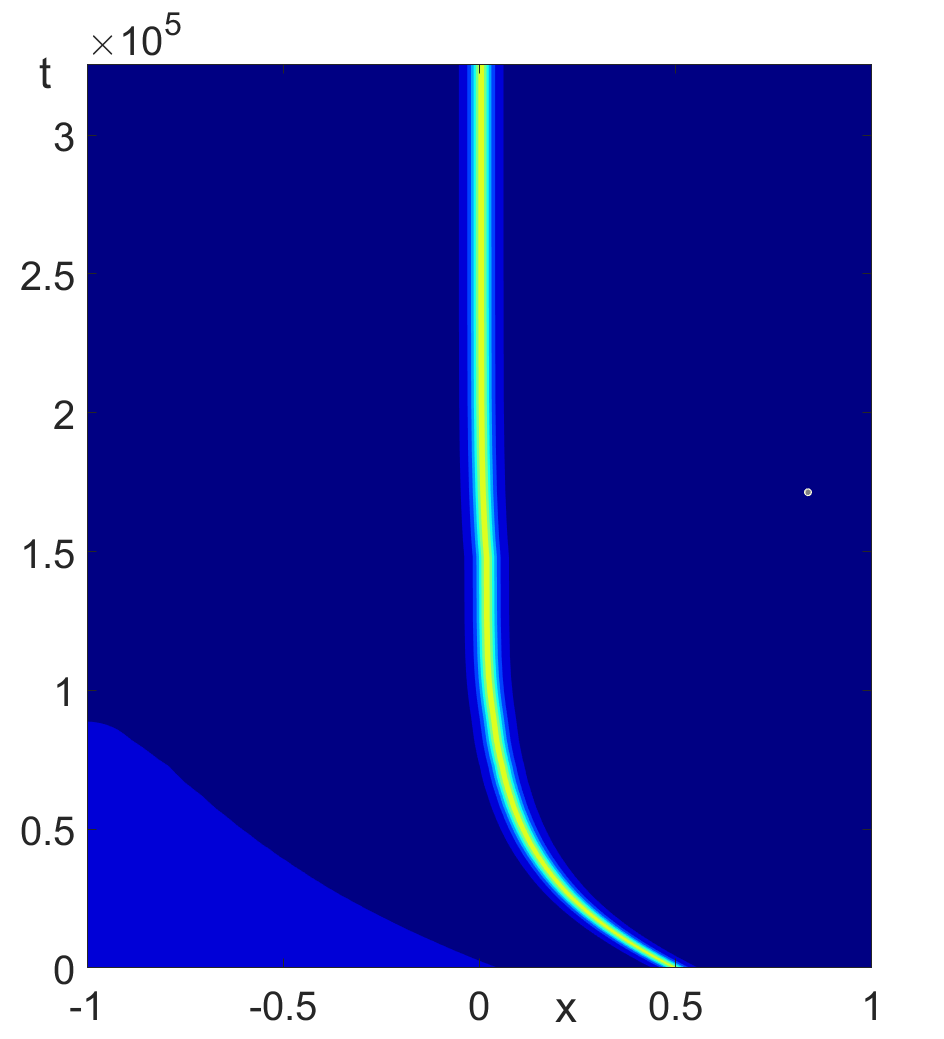}
        \caption{Moving spike when the initial concentration is off-center ($\tau=0.1$).}\label{fig:intro(c)}
    \end{subfigure}
    \caption{Illustration of single spike dynamics in the spatial Solow model \eqref{SS01}. Parameters: $\varepsilon=1\times 10^{-2}, a=1, b=1.$ }
    \vspace{-5mm}
    \label{fig:intro}
\end{figure}

System \eqref{tt1} was originally proposed by \cite{juchem2015capital} as an extension of the spatial Solow model developed in \cite{camacho2004spatial}, which, builds upon the classical Solow–Swan framework \cite{solow1956contribution,swan1956economic}.  The extension introduces labor population dynamics as an additional driver of spatial heterogeneity of the labor density.  This paper extends the analytical framework of \cite{juchem2015capital,kong5177754boundary} by investigating the existence and stability of large-amplitude stationary solutions to system \eqref{tt1} under the small capital diffusion rate regime, see Figure~\ref{fig:intro(a)}.  Such solutions exhibit localized concentration patterns and rich spatial dynamics, including periodic oscillations (Figure~\ref{fig:intro(b)}) and the slow drift motion (Figure~\ref{fig:intro(c)}).  These results and system (\ref{tt1}) offer a rigorous mathematical framework for interpreting the economic agglomeration phenomena, whch is a fundamental feature of spatial economics—the spatial clustering of firms and populations that collectively contribute to higher regional productivity, etc.   A wide range of mechanisms have been investigated to account for the endogenous formation of spatial structures, see \cite{fujita1996role,glaeser2004cities,moses1958location,venables1996equilibrium,weber1929alfred,hirschman1958strategy,perroux1950economic}.  Building on the insights of Marshall \cite{marshall1890principles}, Krugman \cite{krugman1991increasing} first provided a comprehensive general equilibrium framework showing that the interaction among increasing returns, transportation costs, and labor mobility can generate the economic agglomerations and core–periphery patterns.  Further theoretical and empirical perspectives on the mechanisms shaping economic geography and agglomeration can be found in \cite{evans2003development,fingleton2010neoclassical,garretsen2010rethinking,mccann2019theories}.

\subsection{Main results}  This paper extends the results of \cite{kong5177754boundary} and investigates the formation of stable interior spikes in the regime of sufficiently small capital diffusivity.  To perform our analysis, we introduce the dimensionless variables:
\begin{equation}
\begin{aligned}
    \tilde{l}:&= \frac{\chi  }{ d_l }  \left( \frac{1-c}{\delta}  \right) ^{\frac{1}{(1-\theta)}}l, \quad  
\tilde{k}:=\frac{\chi}{d_l} k, \quad  
\tilde{x}:=  \frac{x}{\ell _b} ,\quad 
\tilde{t}=\delta t,\quad \tau:=\tilde \ell_b^2, \\
\tilde \ell_b:&=\sqrt \frac{\delta}{d_l} \ell_b ,\quad 
\tilde a :=\frac{a \chi}{bd_l} \left( \frac{1-c}{\delta}  \right) ^{\frac{1}{(1-\theta)}},\quad
\tilde b: = \frac{b d_l\tilde \ell_b^2}{\chi \delta \left( \frac{1-c}{\delta}  \right) ^{\frac{1}{(1-\theta)}}},
\quad  \varepsilon^2:= \frac{d_k}{d_l\tilde \ell_b^2}, 
\end{aligned}
\end{equation}
which transform  \eqref{tt1} into the following normalized problem:
\begin{equation} \label{eq0}
\left\{
\begin{array}{ll}
\tau \tilde l_{\tilde t} = \tilde l_{\tilde x\tilde x} -  (\tilde l\tilde k_{\tilde x})_{\tilde x} + \tilde b \tilde l({\tilde a} - \tilde l), & \tilde x \in (-1 ,1), \tilde t>0, \\
\tilde k_{\tilde t} = \varepsilon^2 \tilde k_{\tilde x\tilde x} - \tilde k + \tilde k^{\theta}\tilde l^{1-\theta}, & \tilde x \in (-1,1), \tilde t>0,\\
\tilde l_{\tilde x} = \tilde k_{\tilde x} = 0, & \tilde x = -1,1, ~\tilde t>0,\\
\tilde l(\tilde x,0)=\tilde l_0(\tilde x)\geq \not\equiv 0,  \tilde k(\tilde x,0)=\tilde k_0(\tilde x)\geq\not\equiv 0,&\tilde x \in (-1,1).
\end{array}
\right.
\end{equation} 
Dropping the tildes for notational simplicity, one gets the following non-dimensionalized system:
\begin{align}\label{SS01}
\left\{
\begin{array}{ll}
\tau l_t = l_{xx} -  (l k_x)_x + b l({a} - l), & x \in (-1,1), \\
k_t = \varepsilon^2 k_{xx} - k + k^{\theta}l^{1-\theta}, & x \in (-1,1), \\
l_x(\pm 1) = k_x(\pm 1) = 0.
\end{array}
\right.
\end{align}

Concerning $\varepsilon^2 \ll 1$, we study the existence, stability and slow drift dynamics of spiky solutions to (\ref{SS01}).  To begin with, we state the first set of our main results below.
\begin{proposition}\label{prop:1}
Suppose $ab < \tfrac{\pi^2}{4}$. Then there exists a steady state of \eqref{SS01} on $[-1,1]$ that is even with respect to the origin and given by
\begin{equation}\label{prop1asymmainresult}
    l(x) \sim 
    \begin{cases}
       L_0(x/\varepsilon), 
        & |x| \leq \mathcal{O}(\varepsilon), \\[6pt]
        \dfrac{S \cos(\sqrt{ab}\,(|x|-1))}{\cos(\sqrt{ab})}, 
        & \varepsilon \ll |x| \leq 1,
    \end{cases}\quad
    k(x) \sim 
    \begin{cases}
        K_0(x/\varepsilon),
        &  |x| \leq \mathcal{O}(\varepsilon), \\[6pt]
        \dfrac{S \cos(\sqrt{ab}\,(|x|-1))}{\cos(\sqrt{ab})}, 
        & \varepsilon \ll |x| \leq 1,
    \end{cases}
\end{equation}
where $L_0(x/\varepsilon) := S e^{-S} e^{K_0(x/\varepsilon)}$, $S$ is determined by 
\begin{equation}\label{fullinnermatching11}
    S \sqrt{ab} \tan(\sqrt{ab})
    \sim \varepsilon b 
   ~\mathcal{I}(S),\quad \text{with} \quad \mathcal{I}(S):=-\int_{0}^\infty \big[a L_0 - L_0^2 - (aS-S^2)\big]\,dy,
\end{equation} 
and $K_0$ is the solution of
\begin{equation}\label{K0eqprop11}
0 = K_{0yy} - K_0 + S^{1-\theta} K_0^{\theta} e^{(1-\theta)(K_0-S)},
\qquad 
K_{0y}(0)=0,\quad K_{0y}(\infty)=0.
\end{equation}
\end{proposition}

The second set of results addresses the stability of the single spike steady state given by \eqref{prop1asymmainresult}, and is stated as follows.
\begin{proposition}\label{prop:2mainresult}
For the single interior spike solution from \eqref{prop1asymmainresult},  the spike is linearly stable if $\tau = 0$, and undergoes a Hopf bifurcation when $\tau$ becomes sufficiently large.
\end{proposition}
Propositions \ref{prop:1} and \ref{prop:2mainresult} demonstrate that (\ref{SS01}) admits a stable single interior spike in the asymptoic limit of small captial diffusion.  Moreover, we investigate the slow drift dynamics of the quasi-equilibrium counterpart of (\ref{prop1asymmainresult}) and obtain the last set of results, which can be stated as follows 

\begin{proposition}\label{propositon3mainresults}
  For (\ref{SS01}), assume
  that the quasi-equilibrium pattern is linearly stable with respect
  to the large eigenvalues. Then, the
  slow dynamics of the location $x_0$ of the single spike satisfies
  the DAE system:
\begin{align*}
  \left\{\begin{array}{ll}
  \frac{dx_0}{dt}\sim \varepsilon^2S_0\sqrt{ab}\frac{    {(1-\theta)} \int_0^\infty \frac{1}{L_0(z)} \left( K_{0y}^2(z) - (K_0^2(z) - S_0^2) \right)~dz}{2\int_{-\infty}^{\infty} K_{0y}^2 dy}
  \left(\tan[\sqrt{ab}(x_0 + 1)] + \tan[\sqrt{ab}(x_0 - 1)] \right)\,\\   2\varepsilon b\mathcal{I}(S_0)
   = - S_0 \sqrt{ab} \left( \tan[\sqrt{ab}(x_0 - 1)] - \tan[\sqrt{ab}(x_0 + 1)] \right),
   \end{array}
   \right.
\end{align*}
where  $L_0(z):=S_0 e^{-S_0} e^{K_0(z)}$, $\mathcal{I}(S)$ is defined in \eqref{fullinnermatching11} and $K_0$ is defined in (\ref{K0eqprop11}) with $S$ replaced by $S_0$.  In particular, the location
$\bar x^0$ of the single spike true steady-state
solution, is the equilibrium point of the slow dynamics and satisfies
\begin{align}\label{ss:balance}
\tan[\sqrt{ab}(\bar x_0 - 1)] - \tan[\sqrt{ab}(\bar x_0 + 1)] =0.
\end{align}
\end{proposition}

 The paper is organized as follows: In Section~\ref{sec2}, we employ the matched asymptotic method to construct the single interior spike steady state of \eqref{SS1} near $x = 0$ for $\varepsilon \ll 1$.   Section~\ref{sec3} is devoted to the analysis of the stability of the single spike we have constructed. Two types of eigenvalues are discussed: the large eigenvalues corresponding to a order $\mathcal{O}(1)$ dynamics and the small eigenvalues associated with the translational instability. Section~\ref{sec4} investigates the slow dynamics of a single spike quasi-equilibrium. We derived a reduced differential-algebraic equation (DAE) system to describe the slow motion of the location of the spike quasi-equilibrium. All theoretical predictions are verified through numerical simulations using the finite element software FlexPDE7 \cite{Flexpde} and Matlab.

\section{The construction of spike solutions}\label{sec2}
We first employ the method of matched asymptotic expansions to construct a steady-state solution of \eqref{SS1} in the form of a symmetric spike centered at the origin.  Equivalently, this can be achieved by constructing a half-spike on the interval $[0,1]$ subject to Neumann boundary conditions, with the corresponding interior spike obtained by even reflection across the origin.
At equilibrium, the half-spike steady state satisfies
\begin{align}\label{SS1}
\left\{
\begin{array}{ll}
0 = l_{xx} -  (l k_x)_x + b l({a} - l), & x \in (0,1), \\
0 = \varepsilon^2 k_{xx} - k + k^{\theta}l^{1-\theta}, & x \in (0,1)\end{array}
\right.
\end{align}
with Neumann boundary conditions at $x=0$ and $x=1$.  The spike is localized at $x=0$ with spatial extent of order $\mathcal{O}(\varepsilon)$.  To find the spike profile, we introduce the rescaled coordinate $y = x/\varepsilon$, under which the system \eqref{SS1} transforms into
\begin{align}\label{SSinner}
\begin{cases}
0 = L_{yy} - (L K_y)_y + \varepsilon^2 b L(a - L), & y \in \mathbb{R}, \\[4pt]
0 = K_{yy} - K + K^{\theta} L^{1-\theta}, & y \in \mathbb{R}, \\[4pt]
L_y(\pm \infty) = K_y(\pm \infty) = 0.
\end{cases}
\end{align}

We expand the solution in powers of $\varepsilon^2$ as
\begin{equation}\label{ssexpand}
L = L_0 + \varepsilon^2 L_1 + \cdots, 
\qquad 
K = K_0 + \varepsilon^2 K_1 + \cdots.
\end{equation}
Substituting  \eqref{ssexpand} into \eqref{SSinner} and collecting $\mathcal{O}(1)$ terms yield
\begin{equation}\label{inner0}
\begin{cases}
0 = L_{0yy} - (L_0 K_{0y})_y, & y \in (0,\infty), \\[4pt]
0 = K_{0yy} - K_0 + K_0^{\theta} L_0^{1-\theta}, & y \in (0,\infty), \\[4pt]
L_{0y}(0) = K_{0y}(0) = L_{0y}(\infty) = K_{0y}(\infty) = 0.
\end{cases}
\end{equation}
Solving the first equation of \eqref{inner0} with the prescribed far-field behavior $L_0(\pm \infty) = K_0(\pm \infty) = S$, we obtain
\begin{equation}\label{L0}
L_0 = S e^{-S} e^{K_0},
\end{equation}
where $S$ is a constant determined by matching with the outer solution.  Substituting \eqref{L0} into the second equation of \eqref{inner0} yields
\begin{equation}\label{K0}
0 = K_{0yy} - K_0 + S^{1-\theta} K_0^{\theta} e^{(1-\theta)(K_0-S)},
\qquad 
K_{0y}(0)=0,\quad K_{0y}(\infty)=0.
\end{equation}

\begin{figure}
    \centering
    \begin{subfigure}[b]{0.45\linewidth}
        \centering
        \includegraphics[width=\linewidth]{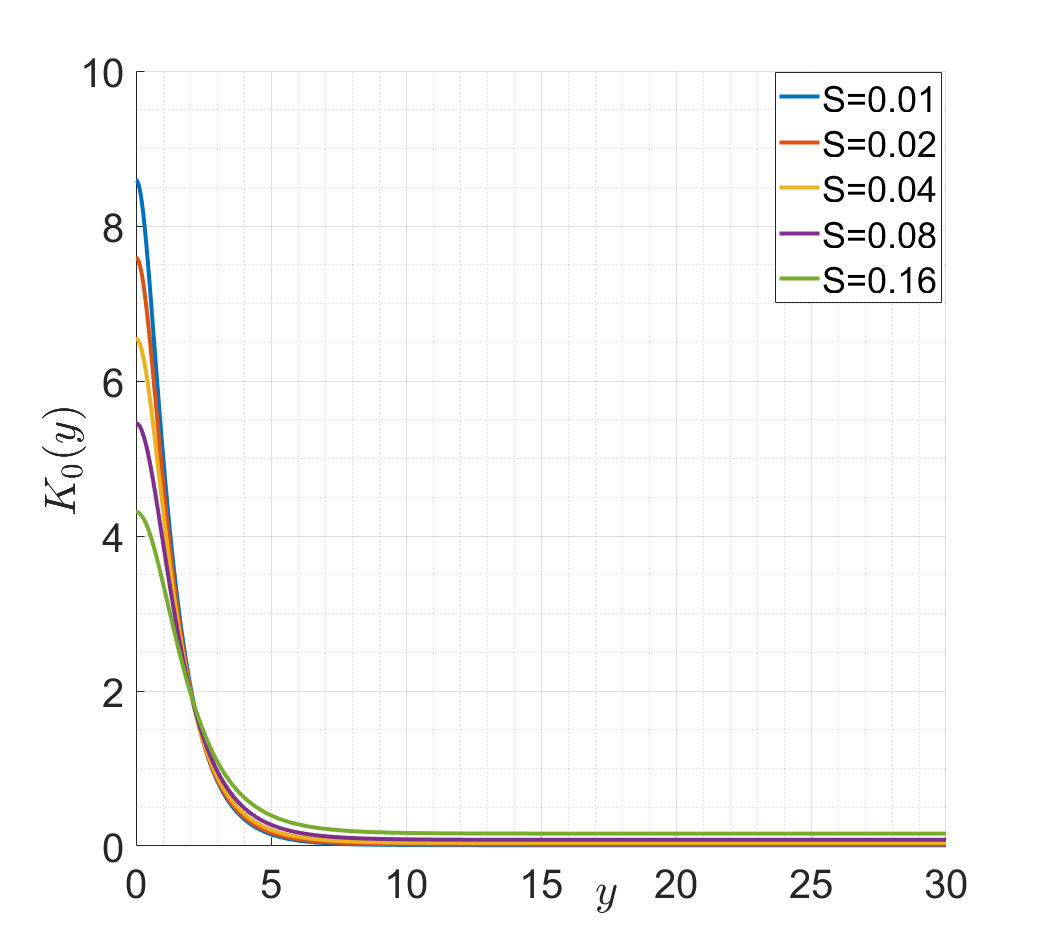}
        \caption{Profile of $K_0(y)$ for various values of $S$ at $\theta=0.5$, computed numerically from \eqref{K0}.}\label{fig:IS(a)}
    \end{subfigure}
    \hfill
    \begin{subfigure}[b]{0.5\linewidth}
        \centering
        \includegraphics[width=\linewidth]{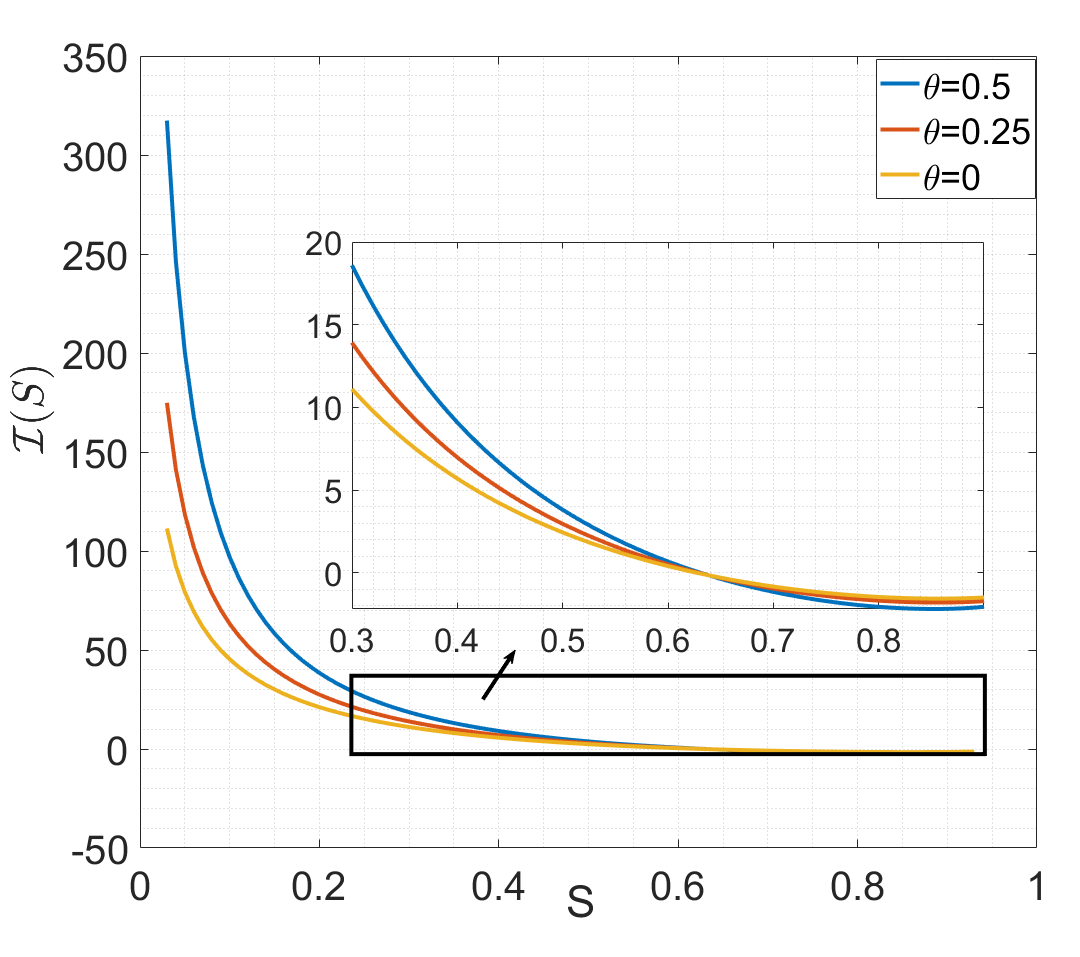}
        \caption{The quantity $\mathcal{I}(S)$ defined in \eqref{IS} for various $\theta$.} \label{fig:IS(b)}
    \end{subfigure}
    \caption{Numerical results illustrating the divergence of $K_0(0)$ and $\mathcal{I}(S)$ as $S\to 0$.}
    \label{fig:IS}
\end{figure}

Equation \eqref{K0} admits two constant solutions, namely $0$ and $S$.  
One can show that there exist a homoclinic solution satisfying $K_0(y) > 0$ and $K_0(\pm \infty) = S$ for $0<S<1$.  
Numerically, \eqref{K0} can be solved as a function of $S$, see Figure~\ref{fig:IS}.  Analytically, it can be expressed in the first integral form
\begin{equation}\label{1int}
K_{0y}^2 - (K_0^2 - S^2) 
+ 2 S^{1-\theta} \int_S^{K_0} z^{\theta} e^{(1-\theta)(z-S)} \, dz = 0.
\end{equation}

In order $\varepsilon^2$, we obtain
\begin{equation}\label{inner1}
\begin{cases}
- b(a L_0 - L_0^2) = L_{1yy} - (L_1 K_{0y} + L_0 K_{1y})_y, & y \in (0,\infty), \\[4pt]
0 = K_{1yy} - K_1 + \theta K_0^{\theta-1} L_0^{1-\theta} K_1 
+ (1-\theta) K_0^{\theta} L_0^{-\theta} L_1, & y \in (0,\infty).
\end{cases}
\end{equation}
From the first equation in \eqref{inner1}, we deduce
\begin{equation}\label{solv20}
 L_{1y} - (L_1 K_{0y} + L_0 K_{1y})
=  \int_0^y \big( a L_0 - L_0^2 \big) \, dy.
\end{equation}
Since $K_{0y} \to 0$, $L_0 \to S$, and $K_{1y} \to L_{1y}$ as $y \to +\infty$, the equation \eqref{solv20} implies
\begin{equation}\label{solv2}
(1-S) L_{1y}(y) \to - b(aS - S^2) y 
- b \int_0^{\infty}  (a L_0 - L_0^2) - (aS - S^2) \, dy,
\qquad y \to \infty.
\end{equation}

In the outer region, sufficiently far from the spike center (i.e., for $x \gg \varepsilon$), we have $k(x) \sim l(x)$  according to the second equation of \eqref{SS1}.  Then the function $l(x)$ satisfies the following boundary-value problem
\begin{equation}\label{outer}
    l_{xx} - (l l_x)_x + b l(a-l) \sim 0,
    \qquad l_x(1) = 0, 
    \qquad l(0) = S.
\end{equation}
For general $S$, the problem \eqref{outer} cannot be solved analytically and must instead be treated numerically.  
Near $x=0$, the local behavior of $l$ is given by
\begin{equation}
    l(x) \sim l(0) + l_x(0^+) x + \tfrac{1}{2} l_{xx}(0^+) x^2 + \cdots.
\end{equation}
In the local coordinate, $x = \varepsilon y$, this expansion becomes
\begin{equation}\label{local}
    l(x) \sim l(0) + l_x(0^+) \varepsilon y + \tfrac{1}{2} l_{xx}(0^+) \varepsilon^2 y^2 + \cdots.
\end{equation}
Meanwhile, the far-field behavior of the inner solution is
\begin{equation}\label{innerfar}
L_y(y) = S 
- \varepsilon^2 y \, \frac{b \int_{0}^\infty \big(a L_0 - L_0^2 - (aS-S^2)\big)\,dy}{1-S} 
- \varepsilon^2 \frac{b(aS-S^2) y^2}{2(1-S)} + \cdots,
\qquad y \to \infty.
\end{equation}
Matching the local outer expansion \eqref{local} with the inner far-field expansion \eqref{innerfar} yields the relation
\begin{equation}\label{matching}
    l_x(0^+) = - \frac{\varepsilon b}{1-S} \int_{0}^\infty \big(a L_0 - L_0^2 - (aS-S^2)\big)\,dy.
\end{equation}
Since the outer solution $l(x)$ depends on $S$ through the boundary condition $l(0) = S$, the matching condition \eqref{matching} provides an implicit equation for $S$, which must be solved numerically.  

In the singular limit $\varepsilon \to 0$, the matching condition \eqref{matching} indicates that $S \to 0$.  
To see this, suppose $S \sim \mathcal{O}(1)$. Then the left-hand side of \eqref{matching} is $\mathcal{O}(1)$, whereas the right-hand side is only $\mathcal{O}(\varepsilon)$, resulting in no consistent solution.  
Thus, we assume $S \ll 1$, which implies $l(x) \ll 1$.  With this assumption, the outer problem can be approximated by the linear equation
\begin{equation}\label{outer_appro}
    l_{xx} + ab\, l \sim 0,
    \qquad l_x(1) = 0, 
    \qquad l(0) = S.
\end{equation}
The solution of \eqref{outer_appro} is
\begin{equation}
    l(x) \sim \frac{S \cos\!\big(\sqrt{ab}\,(x-1)\big)}{\cos\!\big(\sqrt{ab}\big)},
    \qquad 0 < x < 1.
\end{equation}
We remark that, in order to ensure the positivity of $l(x)$ in the outer region, the condition $ab<\frac{\pi^2}{4}$ is required, as stated in Proposition~\ref{prop:1}.

Thus, the dominant part of the matching condition \eqref{matching} becomes
\begin{equation}\label{fullinnermatching}
    S \sqrt{ab} \tan(\sqrt{ab})
    \sim \varepsilon b 
   ~ \mathcal{I}(S),
\end{equation}
where $\mathcal{I}(S)$ is a function of $S$, defined as
\begin{equation}\label{IS}
    \mathcal{I}(S):=-\int_{0}^\infty \big(a L_0 - L_0^2 - (aS-S^2)\big)\,dy.
\end{equation}
For a general $S$, we have to evaluate $\mathcal{I}(S)$ numerically. Figure.~(\ref{fig:IS(a)}) shows the dependence of $\mathcal{I}(S)$ on $S$. We can further estimate $\mathcal{I}(S)$ in \eqref{fullinnermatching} in the limit $S \to 0$.  

Numerical evidence, Figure~\ref{fig:IS(b)}, shows that $K(0) \to \infty$ as $S \to 0$.  
This motivates us to seek the leading-order approximation of \eqref{1int} at $y=0$.  
Let $\xi := K(0)$. We assume $\xi \gg 1$ and $S \ll 1$.  
Evaluating \eqref{1int} at $y=0$ gives
\begin{equation}\label{intappro0}
    \xi^2 - S^2 
    + 2S^{1-\theta} \int_S^\xi z^\theta e^{(1-\theta)(z-S)}\,dz = 0,
    \qquad 0 < \theta < 1.
\end{equation}
Integrating by parts and retaining the dominant contribution, we obtain
\begin{equation}\label{intappro}
    -\frac{1}{2} \xi^2 + \frac{1}{1-\theta} S^{1-\theta} \xi^\theta e^{(1-\theta)\xi} 
    \Bigg(1 + \mathcal{O}\!\left(\frac{1}{(1-\theta)\xi}\right)\Bigg) \sim 0,
    \qquad 0 < \theta < 1.
\end{equation}

\begin{figure}
    \centering
    \includegraphics[width=\linewidth]{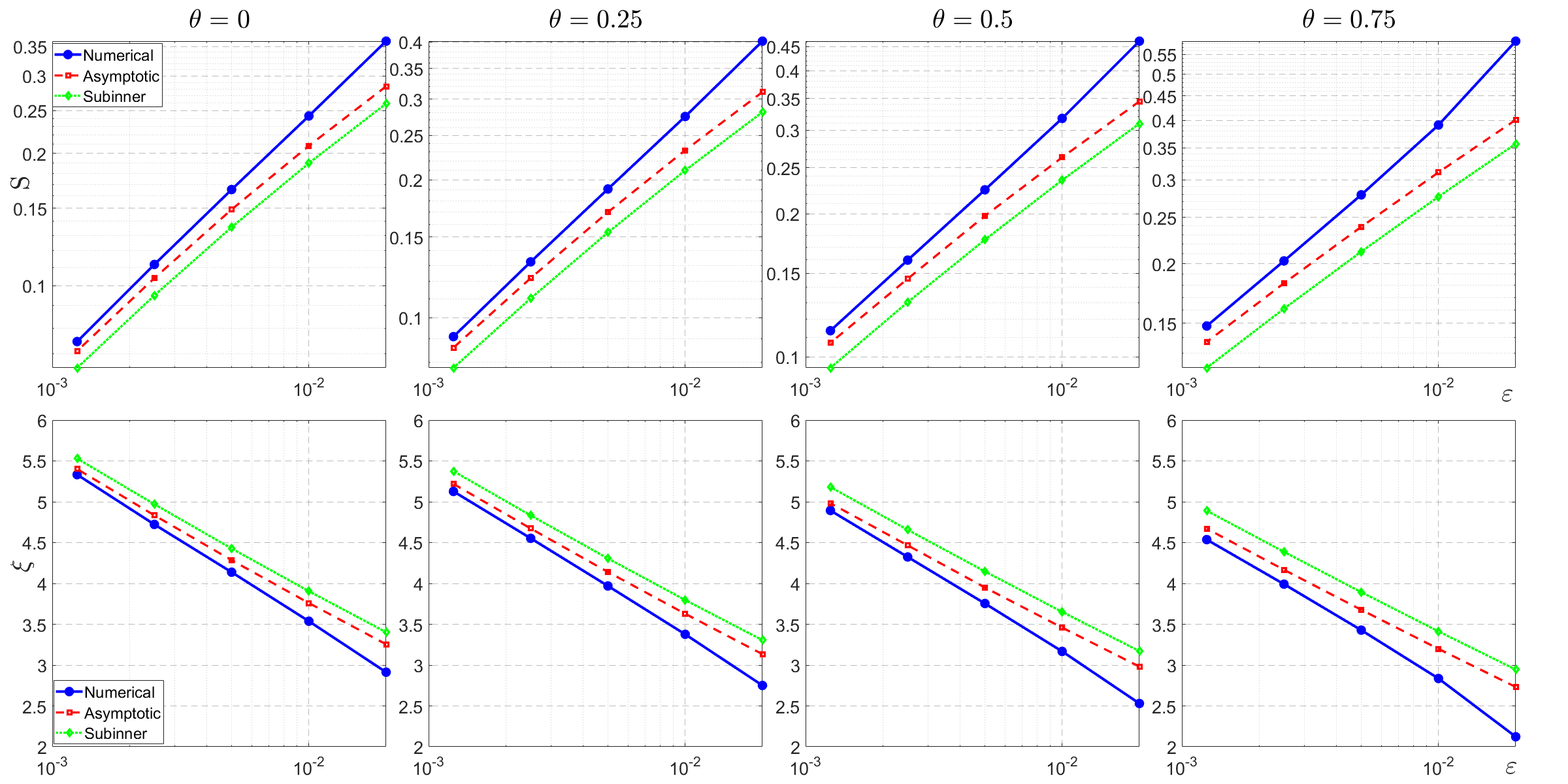}
    \caption{Comparison of $S$ (top panels) and $\xi$ (bottom panels) obtained from different approaches for 
    $\varepsilon = 2\times10^{-2},\; 1\times10^{-2},\; 5\times10^{-3},\; 2.5\times10^{-3},\; 1.25\times10^{-3}$, 
    with parameters $a=1,\; b=0.25$. 
    The $x$-axis is plotted on a logarithmic scale; $S$ is also on a logarithmic scale, while $\xi$ is on a linear scale. 
    Blue curves correspond to direct numerical solutions of the full PDE system \eqref{SS01} using FlexPDE, 
    red curves are obtained from the numerical solution of \eqref{fullinnermatching}, 
    and green curves from the solution of \eqref{subalge}. 
    The agreement between all approaches improves as $\varepsilon$ decreases.}
    \label{fig:Epsi}
\end{figure}

To estimate $\mathcal{I}(S)$ in \eqref{fullinnermatching} in the limit $\xi \gg 1$, we expand $K_0$ in the ``sub-inner'' region: $y \sim \mathcal{O}(\xi^{-1}) \ll 1$.  
Introducing the rescaled variable $z = \tfrac{(1-\theta)\xi y}{2}$, we expand
\begin{equation}\label{K0expand}
    K_0(y) = \xi + \frac{1}{1-\theta} K_{00}(z) + \mathcal{O}(\xi^{-1}).
\end{equation}
Substituting \eqref{K0expand} into \eqref{K0}, at leading order we obtain
\begin{equation}\label{K_00}
    \frac{1-\theta}{4} \,\xi^2 K_{00zz}(z) 
    - \xi^\theta S^{1-\theta} e^{-(1-\theta)S} e^{(1-\theta)\xi} e^{K_{00}} 
    \sim 0.
\end{equation}
From \eqref{intappro}, we obtain
\begin{equation}\label{xi2theta}
    \xi^\theta S^{1-\theta} e^{-(1-\theta)S} e^{(1-\theta)\xi}
    \Big(1 + \mathcal{O}(\xi^{-1})\Big) \sim \frac{(1-\theta)\xi^2}{2}.
\end{equation}
Substituting \eqref{xi2theta} into \eqref{K_00} and neglecting order of $\mathcal{O}(\xi^{-1})$ yields
\begin{equation}\label{K00}
    K_{00zz} + 2 e^{K_{00}} = 0, 
    \qquad K_{00}(0) = 0.
\end{equation}
The solution to \eqref{K00} is
\begin{equation}
    K_{00}(z) = \log\!\left(\sech^2(z)\right).
\end{equation}
It follows that
\begin{equation}
    L_{00}(z) \sim S e^{-S} e^{\xi} \sech^{\tfrac{2}{1-\theta}}(z).
\end{equation}
Therefore, we estimate
\begin{equation}\label{leadint}
\begin{aligned}
    \int_{0}^\infty L_0\,dy 
    &\sim \frac{2 S e^{-S} e^{\xi}}{(1-\theta)\xi}
    \int_{0}^\infty \sech^{\tfrac{2}{1-\theta}}(z)\,dz
    = \frac{2 S e^{-S} e^{\xi}}{(1-\theta)\xi}
    \frac{\sqrt{\pi}\,\Gamma\!\left(1 + \tfrac{1}{1-\theta}\right)}{2 \Gamma\!\left(\tfrac{1}{2} + \tfrac{1}{1-\theta}\right)}, \\[6pt]
    \int_{0}^\infty L_0^2\,dy 
    &\sim \frac{2 S^2 e^{-2S} e^{2\xi}}{(1-\theta)\xi}
    \int_{0}^\infty \sech^{\tfrac{4}{1-\theta}}(z)\,dz
    = \frac{2 S^2 e^{-2S} e^{2\xi}}{(1-\theta)\xi}
    \frac{\sqrt{\pi}\,\Gamma\!\left(1 + \tfrac{2}{1-\theta}\right)}{4 \Gamma\!\left(\tfrac{1}{2} + \tfrac{2}{1-\theta}\right)}.
\end{aligned}
\end{equation}
Using \eqref{leadint}, the leading-order matching condition \eqref{matching} becomes
\begin{equation}\label{one_spike}
\sqrt{ab}\,S\tan(\sqrt{ab})
\sim  \varepsilon b \frac{2 S^2 e^{-2S} e^{2\xi}}{(1-\theta)\xi}
\int_{0}^\infty \sech^{\tfrac{4}{1-\theta}}(z)\,dz.
\end{equation}
Combining \eqref{one_spike} with \eqref{intappro} yields the algebraic system for the $\xi$ and $S$
\begin{equation}\label{subalge}
\begin{aligned}
    \xi^2 &\sim \frac{2}{1-\theta} S^{1-\theta} \xi^\theta e^{(1-\theta)\xi}, \\[6pt]
    S\sqrt{ab}\tan(\sqrt{ab}) 
   &  \sim \varepsilon b \inti{L_0^2}\sim  b\varepsilon\frac{2 S^2 e^{-2S} e^{2\xi}}{(1-\theta)\xi} 
    \int_{0}^\infty \sech^{\tfrac{4}{1-\theta}}(z)\,dz.
\end{aligned}
\end{equation}
Solving \eqref{subalge} asymptotically, we obtain the asymptotic estimates
\begin{equation}\label{asyms}
    \xi \sim \mathcal{O}(|\ln \varepsilon|) \gg 1,
\qquad 
\varepsilon \ll S \sim \mathcal{O}\!\left(\varepsilon |\ln \varepsilon|^{\tfrac{2-\theta}{1-\theta}}\right) \ll 1,
\qquad \varepsilon \to 0,
\end{equation}
which are consistent with our assumptions.  Note that the higher-order terms neglected in \eqref{subalge} are of the order $\mathcal{O}((1-\theta)^{-1}\xi^{-1})$, so the asymptotics are valid provided $(1-\theta)|\ln \varepsilon| \gg 1$.  We summarize the construction of a single spike solution as Proposition \ref{prop:1}.


\begin{figure}
    \centering
    \includegraphics[width=\linewidth]{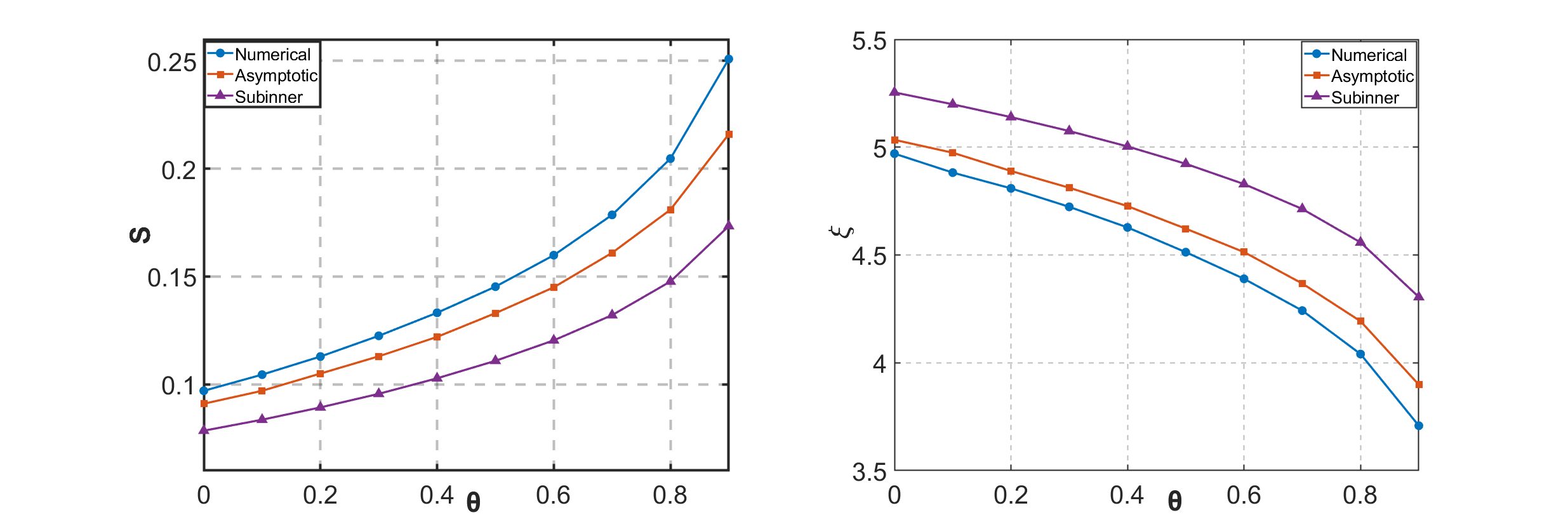}
    \caption{Comparison of $S$ and $\xi$ obtained from different approaches for varying values of $\theta$, with parameters $a=1,\; b=1,\; \varepsilon=2.5\times10^{-3}$. 
    All methods show consistent qualitative trends as $\theta$ increases. 
    For fixed $\varepsilon$, the accuracy of the asymptotic approximation decreases with larger $\theta$, in agreement with the scaling \eqref{asyms}, which indicates that neglected terms are of order $S^2 \sim \varepsilon^2 |\ln \varepsilon|^{\frac{4-2\theta}{1-\theta}}$.}
    \label{fig:theta}
\end{figure}

We validate our asymptotic predictions by comparing the expression \eqref{prop1asymmainresult} with the numerical solution of \eqref{SS1}. 
Figure \ref{fig:Epsi} presents the dependence of $S$ and $\xi$ on $\varepsilon$, showing results from full system simulations, the asymptotic approximation, and the sub-inner approximation. In addition, Figure \ref{fig:theta} illustrates how $S$ and $\xi$ 
vary with $\theta$. The discrepancy between the full numerical results and the asymptotic prediction increases 
as $\theta \to 1$, since the higher-order terms neglected in the expansion are of order $\mathcal{O}(S^2)\sim \mathcal{O} \left(\varepsilon^2 |\ln \varepsilon|^{\frac{4-2\theta}{1-\theta}}\right) $, which becomes more  significant in this limit.

\section{The stability of a interior spike solution}\label{sec3}
We now investigate the stability of the spike equilibrium solution constructed in Proposition~\ref{prop:1}. To analyze stability, we introduce small perturbations of the form
\begin{equation}
    l = l_{s} + e^{\lambda t}\phi(x), \qquad k = k_{s} + e^{\lambda t}\psi(x), \label{perb}
\end{equation}
where $\lambda$ is the spectral parameter and $(\phi,\psi)$ are the corresponding eigenfunctions. Substituting \eqref{perb} into system \eqref{eq0} yields the eigenvalue problem
\begin{subequations}\label{eigen}
\begin{align}
    \tau \lambda \phi &= \phi_{xx} - (\phi k_{sx} + l_s \psi_x )_x + b(a\phi - 2l_s\phi), \label{phi}\\
    \lambda \psi &= \varepsilon^2 \psi_{xx} - \psi + \theta k_s^{\theta-1} l_s^{1-\theta}\psi + (1-\theta) k_s^\theta l_s^{-\theta}\phi, \label{psi}
\end{align}
\end{subequations}
subject to Neumann boundary conditions.

In the inner region, we rescale $x=\varepsilon y$ and let
\begin{equation}
    \phi(x) = \Phi(y), \qquad \psi(x) = \Psi(y),
\end{equation}
which transforms \eqref{eigen} into
\begin{subequations}\label{eigen_in_l}
\begin{align}
    \tau \varepsilon^2 \lambda \Phi &= \Phi_{yy} - (\Phi K_{sy} + L_s \Psi_y )_y + \varepsilon^2 b(a\Phi - 2L_s \Phi), \label{phi_in_l}\\
    \lambda \Psi &= \Psi_{yy} - \Psi + \theta K_s^{\theta-1} L_s^{1-\theta} \Psi + (1-\theta) K_s^\theta L_s^{-\theta}\Phi. \label{psi_in}
\end{align}
\end{subequations}

We expand
\begin{equation}
    \Phi(y) = \Phi_0 + \varepsilon^2 \Phi_1 + \cdots, \qquad
    \Psi(y) = \Psi_0 + \varepsilon^2 \Psi_1 + \cdots, \qquad
    \lambda = \lambda_0 + \varepsilon^2 \lambda_1 + \cdots.
\end{equation}
At leading order, we obtain
\begin{subequations}\label{eigen_in0_l}
\begin{align}
    0 &= \Phi_{0yy} - (\Phi_0 K_{0y} + L_0 \Psi_{0y})_y, \label{phi_in0_l}\\
    \lambda_0 \Psi_0 &= \Psi_{0yy} - \Psi_0 + \theta K_0^{\theta-1} L_0^{1-\theta} \Psi_0 + (1-\theta) K_0^\theta L_0^{-\theta}\Phi_0. \label{psi_in0_l}
\end{align}
\end{subequations}

It is important to note that $(L_{0y},K_{0y})$ is always an eigenfunction of \eqref{eigen_in0_l} with eigenvalue $0$. Consequently, the analysis of \eqref{eigen_in0_l} naturally splits into two cases depending on whether $\lambda_0 \neq 0$ (large eigenvalues) or $\lambda_0 = 0$ (small eigenvalues). In what follows, we first address the case of large eigenvalues.

\subsection{Computation of large eigenvalues}

We investigate the large eigenvalue by studying \eqref{eigen_in0_l}. Solving \eqref{phi_in0_l} for $\Phi_0$ yields
\begin{equation}\label{phi0_l}
    \Phi_0(y) = L_0 \Psi_0 + g_0 L_0,
\end{equation}
where $g_{0}$ is a constant determined by the far-field limits $S_\Phi := \Phi_0(+\infty)$ and $S_\Psi := \Psi_0(+\infty)$ via
\begin{equation} \label{g0}
    g_0 = \frac{S_\Phi}{S} - S_\Psi.
\end{equation}
Substituting \eqref{phi0_l} into \eqref{psi_in0_l}, we obtain the following eigenvalue problem:
\begin{equation}\label{Nlep}
    \lambda_0 \Psi_0 = \Psi_{0yy} - \Psi_0 + \theta K_0^{\theta-1} L_0^{1-\theta} \Psi_0 + (1-\theta) K_0^\theta L_0^{1-\theta}\big(\Psi_0 + g_0\big).
\end{equation}

In the outer region, the leading order of \eqref{psi} implies
\begin{equation}
    \psi \sim \frac{1-\theta}{1-\theta + \lambda_0}\phi,
\end{equation}
and hence $\phi$ satisfies
\begin{equation}\label{exactphi}
    \tau \lambda_0 \phi \sim \phi_{xx} - \left(\phi l_{sx} + \frac{1-\theta}{1-\theta+\lambda_0} l_s \phi_x \right)_x + b(a\phi - 2l_s\phi), \qquad \phi'(\pm 1) = 0,
\end{equation}
with the boundary condition at $x=0$ given by $\phi(0) = S_\Phi.$ Since $l_s \ll 1$ and $S \ll 1$,  \eqref{exactphi} reduces to the approximation
\begin{equation}\label{approphi}
    \tau \lambda_0 \phi \sim \phi_{xx} + ab \phi, \qquad \phi(0) = S_\Phi.
\end{equation}
For $ab <\tfrac{\pi^2}{4}$, the nonzero solution of \eqref{approphi} is given by
\begin{equation}\label{phiappro3}
    \phi(x) \sim \frac{S_\Phi}{\cos\mu} \begin{cases}
        \cos\big(\mu(x-1)\big), & 0 < x < 1,\\
        \cos\big(\mu(x+1)\big), & -1 < x < 0,
    \end{cases}\qquad \mu := \sqrt{ab - \tau\lambda_0}.
\end{equation}

\begin{figure}
    \centering
    \begin{subfigure}[t]{0.48\linewidth}
        \centering
        \includegraphics[width=\linewidth]{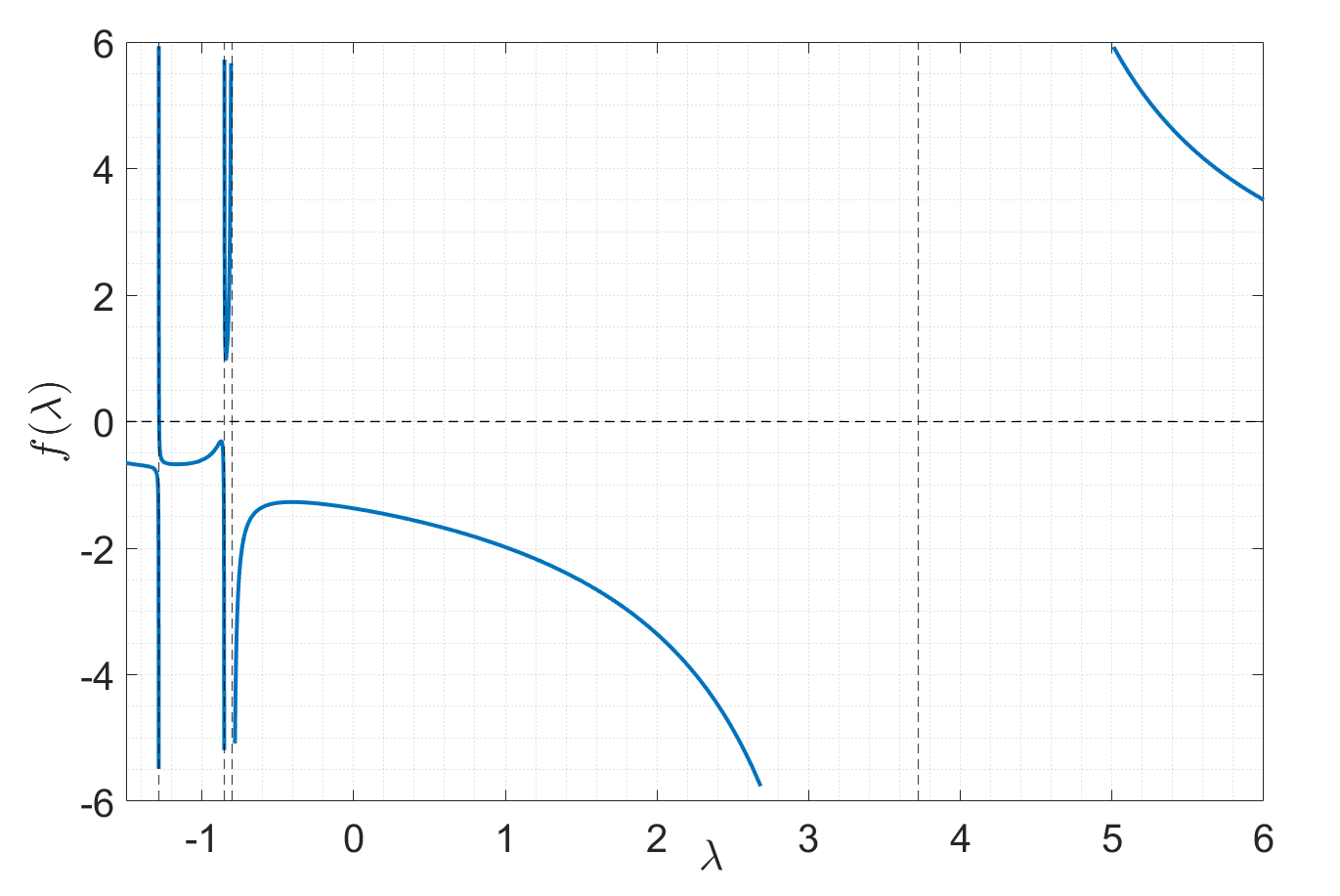}
        \caption{$\theta = 0$.}
        \label{fig:NLEP_f:a}
    \end{subfigure}
    \hfill
    \begin{subfigure}[t]{0.48\linewidth}
        \centering
        \includegraphics[width=\linewidth]{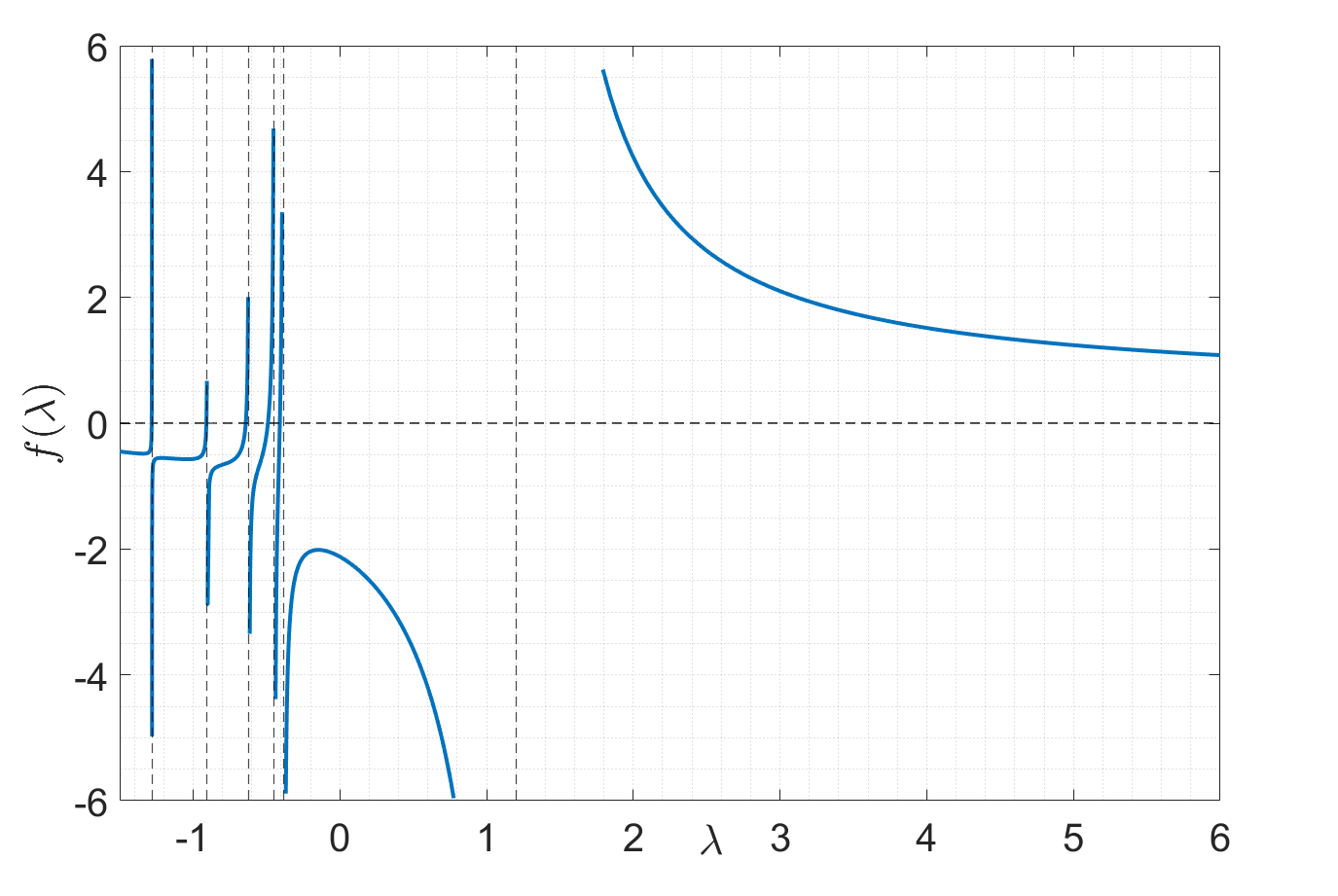}
        \caption{$\theta = 0.5$.}
        \label{fig:NLEP_f:b}
    \end{subfigure}
    \caption{The function $f(\lambda)$ for real $\lambda$ at $\tau=0$, with parameters $a=1,\; b=1,\; \varepsilon=1\times 10^{-2}$. 
    In both cases, no positive real roots of $f(\lambda)=0$ are observed.}
    \label{fig:NLEP_f}
\end{figure}
Matching the local behavior of the outer solution with the far-field behavior of the inner solution, we obtain
\begin{equation}\label{SphiSpsi}
    S_{\Psi}=\psi(0)=\frac{1-\theta }{1-\theta+\lambda_0}\phi(0)= \frac{1-\theta }{1-\theta+\lambda_0} S_\Phi,
\end{equation}
together with the jump condition
\begin{equation*}\label{jumpphi0}
    \left(1- \frac{1-\theta}{1-\theta+\lambda}S\right)\left(\phi_x(0^+)-\phi_x(0^-)\right)
    =S_\Phi\left(l_{sx}(0^+)-l_{sx}(0^-)\right)- \varepsilon b \int_{-\infty}^\infty (a-2L_0)(\Phi_0-S_{\Phi})\,dy.
\end{equation*}
Substituting the approximation \eqref{phiappro3} into the above equation, we obtain
\begin{equation}\label{jumpphi1}
    \left(1- \frac{1-\theta}{1-\theta+\lambda_0}S\right)2S_\Phi\mu \tan \mu
    =2S_\Phi S \sqrt{ab} \tan\sqrt{ab} - \varepsilon b\int_{-\infty}^\infty (a-2L_0)(\Phi_0-S_{\Phi})\,dy,
\end{equation}
Using \eqref{g0} and \eqref{SphiSpsi}, we simplify \eqref{jumpphi1} as
\begin{equation}\label{SPHI}
    2g_0 S\mu \tan \mu= 2S_\Phi S \sqrt{a} \tan\sqrt{a} - \varepsilon b \int_{-\infty}^\infty (a-2L_0)(\Phi_0-S_{\Phi})\,dy.
\end{equation}
Replacing $S_\Phi$ in \eqref{SPHI} with $S(S_\Psi+g_0)$, we obtain
\begin{equation}
    2g_0 S\mu \tan \mu= 2(S_\Psi+g_0) S^2 \sqrt{ab} \tan\sqrt{ab} - \varepsilon b \int_{-\infty}^\infty (a-2L_0)(L_0\Psi_0+g_0L_0-SS_{\Psi}-Sg_0)\,dy.
\end{equation}
Since the term $2(S_\Psi+g_0) S^2 \sqrt{ab} \tan\sqrt{ab}$ is of order $\mathcal{O}(S^2)$, asymptotically small compared to the others, we neglect it and obtain
\begin{equation}\label{g00}
    g_0  = - \frac{\varepsilon b\int_{-\infty}^\infty (a-2L_0)(L_0\Psi_0-SS_{\Psi})\,dy}{2 S\mu \tan \mu +  \varepsilon b\int_{-\infty}^\infty (a-2L_0)(L_0-S)\,dy}.
\end{equation}
Substituting \eqref{g00} into \eqref{Nlep},  we obtain a nonlocal eigenvalue problem
\begin{align}\label{Nlep_1}
    \lambda_0\Psi_0 
= &\Psi_{0yy}-\Psi_0 + \theta K_0^{\theta-1}L_0^{1-\theta} \Psi_0\\
&+(1-\theta) K_0^\theta L_0^{1-\theta}\left(\Psi_0- \frac{\varepsilon b \int_{-\infty}^\infty (a-2L_0)(L_0\Psi_0-SS_{\Psi})\,dy}{2 S\mu \tan \mu + \varepsilon b \int_{-\infty}^\infty (a-2L_0)(L_0-S)\,dy} \right).\nonumber
\end{align}

The eigenvalue problem \eqref{Nlep_1} can be solved numerically. Instead of tackling it directly, we reformulate it as the equivalent scalar condition $f(\lambda_0)=0$, where
\begin{align}
    f(\lambda):=&\varepsilon(1-\theta)b\bigg( \int_{-\infty}^\infty (a-2L_0)(L_0-S) (\mathcal{L}_N-\lambda)^{-1}  (K_0^\theta L_0^{1-\theta})\,dy \\
   & + \int_{-\infty}^\infty (a-2L_0) S \left((\mathcal{L}_N-\lambda)^{-1}(K_0^\theta L_0^{1-\theta}) -S_{\Psi}\right)\,dy \bigg)\nonumber \\
 &  - \left(2 S\mu \tan \mu +  \varepsilon b\int_{-\infty}^\infty (a-2L_0)(L_0-S)\,dy\right)\nonumber,
\end{align}
and $\mathcal{L}_N$ is the operator
\begin{equation}
    \mathcal{L}_N:= \frac{\partial^2 }{\partial y^2} -1 +  \theta K_0^{\theta-1}L_0^{1-\theta} +(1-\theta) K_0^\theta L_0^{1-\theta}.
\end{equation}

\begin{proposition}
    In the limit $\varepsilon \to 0$, the order $\mathcal{O}(1)$ eigenvalue $\lambda$ of the lineraized problem \eqref{eigen_in_l} is asymptotic to the nonlocal eigenvalue problem \eqref{Nlep_1} when $ab < \frac{\pi^2}{4}$. 
\end{proposition}

\begin{figure}[tbp]
    \centering
    \begin{subfigure}[t]{0.35\linewidth}
        \centering
    \includegraphics[width=\linewidth]{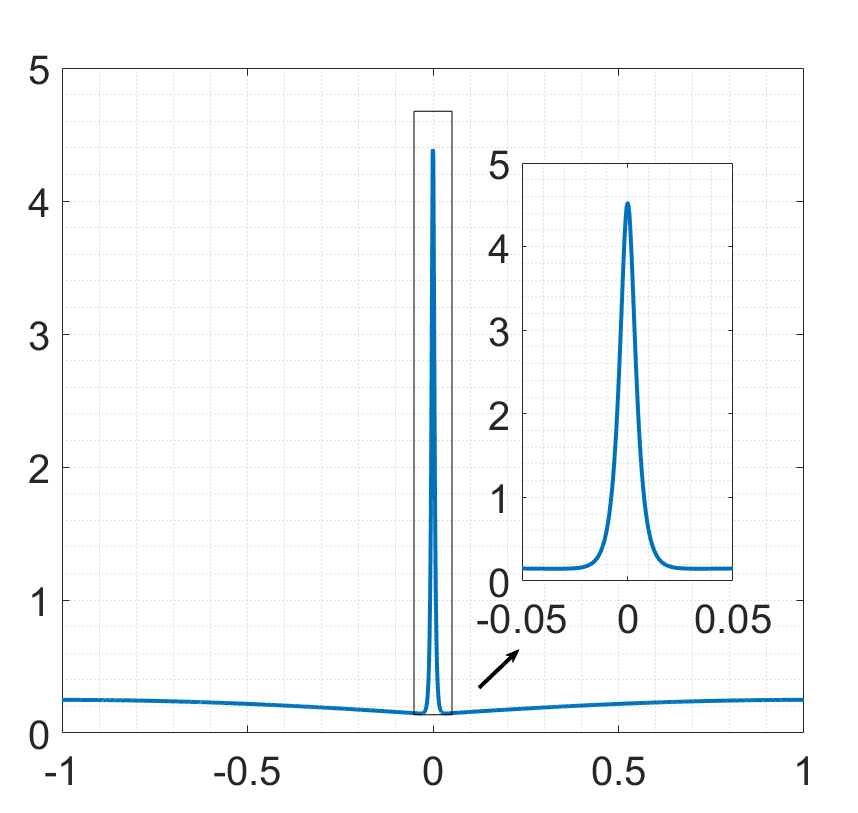}
        \caption{Steady-state profile of $k_s$. }
        \label{fig:Hopfdynamics:a}
    \end{subfigure}
    \hfill
    \begin{subfigure}[t]{0.6\linewidth}
        \centering
    \includegraphics[width=\linewidth]{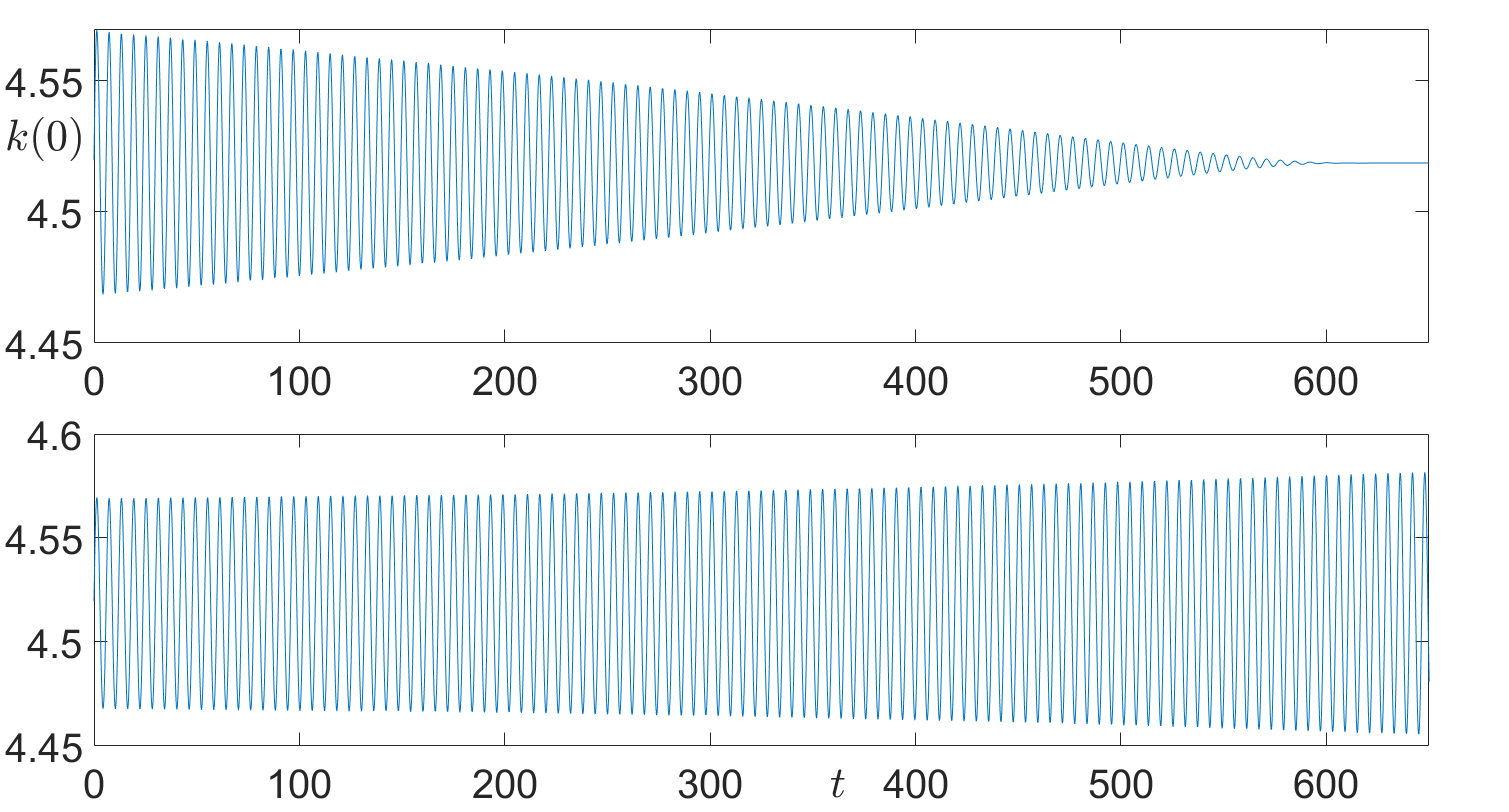}
        \caption{Time evolution of the spike height $k(0)$ before ($\tau=1.43$, top) and after ($\tau=1.44$, bottom) the Hopf bifurcation.}
        \label{fig:Hopfdynamics:b}
    \end{subfigure}
    \caption{Dynamics near the Hopf bifurcation. Paramters are $\varepsilon = 2.5\times10^{-3},\; \theta=0.5,\; a=1,\; b=1$. Simulations are performed using FlexPDE7.}
    \label{fig:Hopfdynamics}
\end{figure}

\begin{figure}
    \centering
    \includegraphics[width=\linewidth]{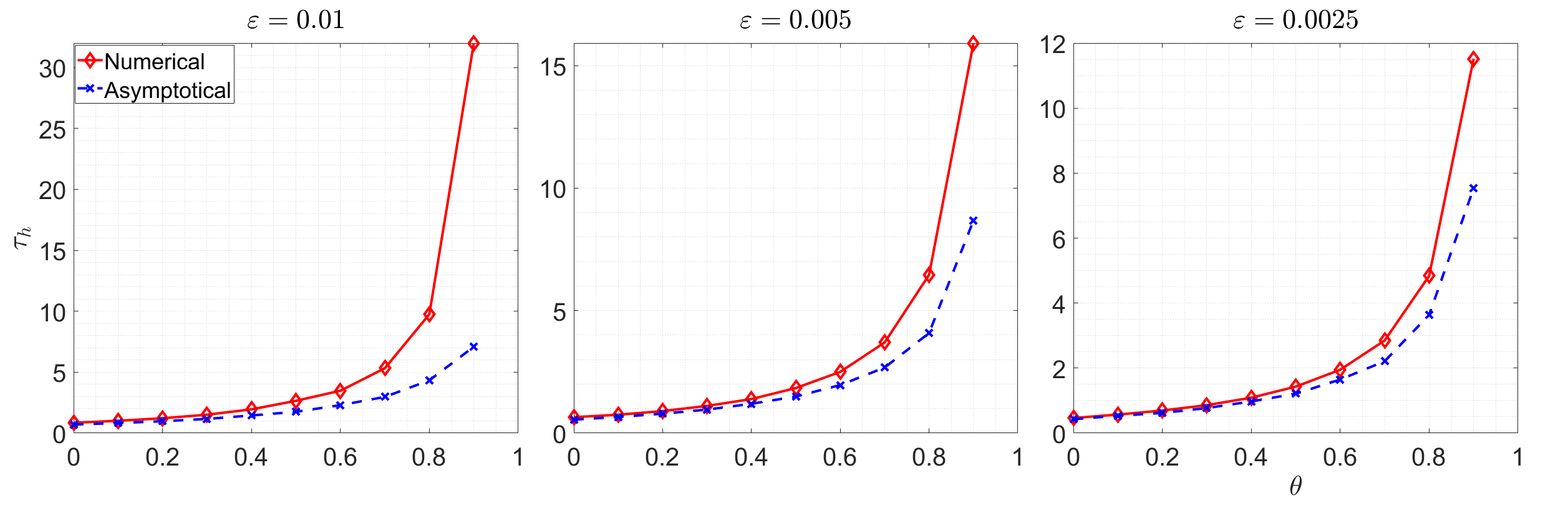}
    \caption{Hopf bifurcation thresholds computed by asymptotical result (NLEP problem \eqref{Nlep_1}) and \eqref{eigen_in_l} at values of $\varepsilon=1\times10^{-2},5\times 10^{-3},2.5\times 10^{-3}$, with parameter $a=1,b=1$.  The agreement between the asymptotic and numerical results improves as $\varepsilon$ decreases and $\theta$ becomes smaller.}
    \label{fig:tauhopf}
\end{figure}

Figure \ref{fig:NLEP_f} shows the profile of $f(\lambda)$ for real values of $\lambda$ at $\tau=0$ for $\theta=0,0.5$, demonstrating that no positive real eigenvalues can be found. Indeed, we can not find a complex eigenvalue with positive real part by numerically solving \eqref{Nlep_1}, indicating that a single interior spike is stable when $\tau=0$, which is in agreement with directly simulating the PDE system \eqref{SS01}.  As $\tau$ increases,  we can find a pure imaginary root to \eqref{Nlep_1} numerically at $\tau_h$, indicating that a Hopf bifurcation occurs at $\tau=\tau_h$.  We numerically compute this $\tau_h$ for various values of $\theta$, and compare it with the eigenvalue computed by directly solving \eqref{eigen} and numerical simulations of full PDE system \eqref{SS01}. Figure~\ref{fig:tauhopf} gives the comparison of $\tau_h$ computed from different approaches at different  $\varepsilon$. The agreement improves as $\varepsilon$ and $\theta$ decrease. 

 In Appendix \ref{app:A}, we further investigate the NLEP problem \eqref{Nlep_1} analytically in the sub-inner region. We show that if a conjecture proposed in \cite{kolokolnikov2009existence} is true, the linearized problem \eqref{eigen_in_l} has no eigenvalue with positive real part at leading order when $\tau=0$.

\subsection{Computation of small eigenvalues}
 
We consider the case of small eigenvalues, scaling them as
\begin{equation}
\lambda = \varepsilon^2 \lambda_1.
\end{equation}

In the inner region, the eigenvalue problem becomes
\begin{subequations}
\label{eigen_in}%
\begin{align}
\tau \varepsilon^4 \lambda_1 \Phi & = \Phi_{yy} - (\Phi K_{sy} + L_s \Psi_y )_y + \varepsilon^2 b (a\Phi L_s - 2L_s\Phi), \label{phi_in} \\
\lambda_1 \varepsilon^2  \Psi  & = \Psi_{yy} - \Psi + \theta K_s^{\theta-1}L_s^{1-\theta} \Psi + (1-\theta) K_s^\theta L_s^{-\theta}\Phi. \label{psi_in_s}
\end{align}
\end{subequations}
At leading order, this system reduces to
\begin{subequations}
\label{eigen_in0}%
\begin{align}
0 & = \Phi_{0yy} - (\Phi_0 K_{0y} + L_0 \Psi_{0y} )_y, \label{phi_in0} \\
0 & = \Psi_{0yy} - \Psi_0 + \theta K_0^{\theta-1}L_0^{1-\theta} \Psi_0 + (1-\theta) K_0^\theta L_0^{-\theta}\Phi_0. \label{psi_in0}
\end{align}
\end{subequations}
A solution to \eqref{eigen_in0} is given by
\begin{equation}\label{phi0sm}
    \Phi_0 = L_{0y}, \quad \Psi_0 = K_{0y}.
\end{equation}

Proceeding to the $\mathcal{O}(\varepsilon^2)$ terms and using \eqref{phi0sm} we obtain the system:
\begin{subequations}
\label{eigen_in2}%
\begin{small}
\begin{align}
-b(aL_{0y} - 2L_0L_{0y}) + (L_{0y} K_{1y} + L_1 K_{0yy} )_y & = \Phi_{1yy} - (\Phi_1 K_{0y} + L_0 \Psi_{1y} )_y, \label{phi_in2} \\
\lambda_1 K_{0y} - \left( \theta K_0^{\theta-1}L_0^{1-\theta} \right)_y L_1 - \left( (1-\theta) K_0^\theta L_0^{-\theta}\right)_y K_1 & = \Psi_{1yy} - \Psi_1 + \theta K_0^{\theta-1}L_0^{1-\theta} \Psi_1 + (1-\theta) K_0^\theta L_0^{-\theta}\Phi_1. \label{psi_in2}
\end{align}
\end{small}
\end{subequations}
The solvability condition of \eqref{eigen_in2} will determine the value of $\lambda_1$.

To facilitate the analysis, we rewrite the system \eqref{eigen_in2} in matrix form using the linear operator $\mathcal{L}$, defined as
\begin{equation}
\mathcal{L}: = \begin{pmatrix}
\displaystyle \frac{\partial^2}{\partial y^2} - K_{0y}\frac{\partial}{\partial y} - K_{0yy} & \displaystyle -L_{0y}\frac{\partial}{\partial y} - L_0 \frac{\partial^2}{\partial y^2} \\[1em]
(1-\theta) K_0^\theta L_0^{-\theta} & \displaystyle \frac{\partial^2}{\partial y^2} - 1 + \theta K_0^{\theta-1}L_0^{1-\theta}
\end{pmatrix}.
\end{equation}
We rewrite \eqref{eigen_in2} as
\begin{align}\label{matrix_raw}
    \mathcal{L} \begin{pmatrix}
        \Phi_1 \\
        \Psi_1
    \end{pmatrix} = \begin{pmatrix}
        0 \\
        \lambda_1 K_{0y}
    \end{pmatrix} + \begin{pmatrix}
        (L_0 K_{1y} + L_1 K_{0y} )_y - b(a L_{0y} - 2L_0L_{0y}) \\[1em]
        \left( \theta K_0^{\theta-1}L_0^{1-\theta} \right)_y L_1 + \left( (1-\theta) K_0^\theta L_0^{-\theta}\right)_y K_1  
    \end{pmatrix}.
\end{align}
We now observe that by differentiating the system \eqref{inner1} with respect to $y$, we find
\begin{equation}\label{L_1y}
\mathcal{L} \begin{pmatrix}
L_{1y} \\
K_{1y}
\end{pmatrix} = \begin{pmatrix}
        (L_0 K_{1y} + L_1 K_{0y} )_y - b(a L_{0y} - 2L_0L_{0y}) \\[1em]
        \left( \theta K_0^{\theta-1}L_0^{1-\theta} \right)_y L_1 + \left( (1-\theta) K_0^\theta L_0^{-\theta}\right)_y K_1  
\end{pmatrix}.
\end{equation}
Using \eqref{L_1y}, we can simplify \eqref{matrix_raw} to
\begin{align}\label{matrix_simple}
    \mathcal{L} \begin{pmatrix}
        \Phi_1 - L_{1y} \\
        \Psi_1 - K_{1y}
    \end{pmatrix} = \begin{pmatrix}
        0 \\
        \lambda_1 \Psi_0
    \end{pmatrix}.
\end{align}
The first equation of \eqref{matrix_simple} can be solved to give
\begin{equation}
    \Phi_1 - L_{1y} = L_0 (\Psi_1 - K_{1y}) + g_1 L_0,
\end{equation}
where $g_1$ is an integration constant. Given the boundary conditions $\Phi_1(\pm\infty) = \Psi_1(\pm\infty) = 0$, we deduce the following far-field conditions:
\begin{equation}\label{far2}
\begin{aligned}
        \Phi_1(\infty) - \Phi_1(-\infty) &= L_{1y}(\infty) - L_{1y}(-\infty), \\
        \Phi_{1y}(\infty) - \Phi_{1y}(-\infty) &= L_{1yy}(\infty) - L_{1yy}(-\infty) = 0.
\end{aligned}
\end{equation}

We now derive the solvability condition for the system \eqref{matrix_simple}. Let $\mathcal{L}^{\dagger}$ denote the adjoint operator of $\mathcal{L}$:
\begin{equation}
\mathcal{L}^\dagger := \begin{pmatrix}
\displaystyle \frac{\partial^2}{\partial y^2} + K_{0y}\frac{\partial}{\partial y} & \displaystyle (1-\theta) K_0^\theta L_0^{-\theta} \\[1em]
\displaystyle -L_{0y}\frac{\partial}{\partial y} - L_0 \frac{\partial^2}{\partial y^2} & \displaystyle \frac{\partial^2}{\partial y^2} - 1 + \theta K_0^{\theta-1}L_0^{1-\theta}
\end{pmatrix}.
\end{equation}
Let $\begin{pmatrix} P \\ Q \end{pmatrix}$ be a vector in the kernel of $\mathcal{L}^\dagger$, satisfying
\begin{equation}\label{PQ}
    \mathcal{L}^{\dagger} \begin{pmatrix} P \\ Q \end{pmatrix} = 0.
\end{equation}
Note that the constant vector $(P_e, Q_e) = (1, 0)$ is a solution to \eqref{PQ}.

We seek a specific solution to \eqref{PQ} subject to the constraints that $P$ and $Q$ are odd functions with $P_y(\infty) = Q_y(\infty) = 0$. The component equations for $P$ and $Q$ are:
\begin{align}
    P_{yy} + K_{0y}P_y + (1-\theta) K_0^\theta L_0^{-\theta} Q &= 0, \label{P} \\
    Q_{yy} - Q + \theta K_0^{\theta-1}L_0^{1-\theta} Q - L_{0y}P_y - L_0 P_{yy} &= 0. \label{Q}
\end{align}
Multiplying \eqref{P} by $L_0$ and adding the result to \eqref{Q} yields a simplified equation for $Q$:
\begin{equation}\label{QQ}
     Q_{yy} - Q + \theta K_0^{\theta-1}L_0^{1-\theta} Q - (1-\theta) K_0^\theta L_0^{1-\theta} Q = 0.
\end{equation}
One solution to \eqref{QQ} is given by
\begin{equation}\label{QQ_S}
    Q = K_{0y}.
\end{equation}
Substituting this into \eqref{P}, we find that $P$ must satisfy
\begin{equation}
L_{0}P_y = - \int_{+\infty}^y (1-\theta) L_0^{1-\theta} K_0^\theta K_{0\xi}  d\xi.
\end{equation}
Consequently, $P$ is given by
\begin{equation}
P(y) = -\int_0^y \frac{1}{L_0(z)} \left( \int_{+\infty}^z (1-\theta) L_0^{1-\theta}(\xi) K_0^\theta(\xi) K_{0\xi}(\xi)  d\xi \right) dz, \quad \text{with} \quad Q(y) = K_{0y}(y).
\end{equation}
To simplify this expression, we use the identity
\begin{equation}\label{P_Indentity}
    K_{0y}^2(z) - (K_0^2(z) - S^2) = -2\int_{\infty}^{z} L_0^{1-\theta}(\xi) K_0^\theta(\xi) K_{0\xi}(\xi)  d\xi,
\end{equation}
which allows us to evaluate $P(y)$ as
\begin{equation}\label{PS}
    P(y) = \frac{(1-\theta)}{2} \int_0^y \frac{1}{L_0(z)} \left( K_{0y}^2(z) - (K_0^2(z) - S^2) \right) dz.
\end{equation}

We now apply the solvability condition to \eqref{matrix_simple} by taking the $\mathbf{L}^2(\mathbb{R})$ inner product of $(P, Q)^T$ with the right-hand side. This yields:
\begin{equation}\label{solva}
    \int_{-\infty}^{\infty} \lambda_1 K_{0y} Q~  dy = \left( P (\Phi_1 - L_{1y})_y - P L_0 (\Psi_1 - K_{1y})_y \right) \Bigg\vert_{-\infty}^{\infty}.
\end{equation}
Note that we have the following far-field behaviors:
\begin{equation} \label{far4}
\begin{aligned}
P(-\infty) &= -P(\infty), \Psi_1(\pm\infty) = \Phi_1(\pm \infty) , \\
L_{0}(\pm \infty) &= S , K_{1yy}(\pm\infty) = L_{1yy}(\pm \infty) = -a b(S - S^2) .
\end{aligned}   
\end{equation}
Using the conditions in \eqref{far4} and neglecting terms of $\mathcal{O}(S^2)$, we evaluate the right-hand side of \eqref{solva}:
\begin{equation}
\left( P (\Phi_1 - L_{1y})_y - P L_0 (\Psi_1 - K_{1y})_y \right) \Bigg\vert_{-\infty}^{\infty} = P(\infty) \left( \Phi_{1y}(+\infty) + \Phi_{1y}(-\infty) + 2abS \right).
\end{equation}
Solving \eqref{solva} for $\lambda_1$ and using $Q = K_{0y}$, we obtain:
\begin{equation}\label{inner_lambda}
 \lambda_1 = \frac{1}{\int_{-\infty}^{\infty} K_{0y}^2 dy } P(\infty) \left( \Phi_{1y}(+\infty) + \Phi_{1y}(-\infty) + 2abS \right).
\end{equation}
The final step is to determine the values of $\Phi_{1y}(\pm \infty)$ by matching with the outer solution.

We now analyze the outer region. Given that $l_s \sim k_s$ in this region, the equation for $\psi$ simplifies to
\begin{equation}
\psi \sim \frac{(1-\theta) k_s^\theta l_s^{-\theta}\phi}{1 - \theta k_s^{\theta-1}l_s^{1-\theta}} = \phi.
\end{equation}
Consequently, the equation for $\phi$ becomes
\begin{equation}
\begin{aligned}
     \phi_{xx} - (\phi l_s)_{xx} + b(a\phi - 2l_s\phi) &\sim 0.
\end{aligned}
\end{equation}
Since $l_s \sim \mathcal{O}(\varepsilon)$, the leading-order equation for $\phi$ is obtained by neglecting the small terms involving $l_s$:
\begin{equation}\label{outer_eig_ss}
          \phi_{xx} + a b \phi \sim 0, \quad \text{with boundary conditions} \quad \phi_x(\pm 1) = 0.
\end{equation}

We match the local behavior of this outer solution $\phi$ near $x=0$ with the far-field behavior of the inner solution given in \eqref{far2}. This matching imposes the following jump conditions at $x=0$:
\begin{equation}
    \begin{aligned}
             \phi(0^+) - \phi(0^-) &\sim \varepsilon^2 \left( L_{1y}(\infty) - L_{1y}(-\infty)\right) \sim \varepsilon^2 2b\mathcal{I}(S) , \\
              \phi_x(0^+) - \phi_x(0^-) &\sim 0.
    \end{aligned} 
\end{equation}
From the steady-state matching condition \eqref{fullinnermatching}, we have the identity
\begin{equation}\label{ssmatching}
     \varepsilon^2 b \mathcal{I}(S)= S\varepsilon \sqrt{ab} \tan(\sqrt{ab}).
\end{equation}
Therefore, the jump in $\phi$ is given by
\begin{equation}\label{eig_jump_c}
     \phi(0^+) - \phi(0^-) \sim 2S\varepsilon \sqrt{ab} \tan(\sqrt{ab}).
\end{equation}

The outer solution $\phi$, satisfying \eqref{outer_eig_ss} and the above jump conditions \eqref{eig_jump_c} at $x=0$, is
\begin{equation}
    \phi(x) = 2\varepsilon S\sqrt{ab}\tan(\sqrt{ab}) \cdot 
    \begin{cases}
         \dfrac{\cos(\sqrt{ab}(x-1))}{2\cos(\sqrt{ab})}, & 0 < x < 1;  \\[1em]
         -\dfrac{\cos(\sqrt{ab}(x+1))}{2\cos(\sqrt{ab})}, & -1 < x < 0.
    \end{cases} 
\end{equation}
From this solution, we compute the required outer derivatives at the origin:
\begin{equation}
    \Phi_{1y}(\pm\infty) = \phi_{x}(0^\pm) = \varepsilon S a b \tan^2(\sqrt{ab}).
\end{equation}
Substituting this result into the expression for the eigenvalue \eqref{inner_lambda} yields
\begin{equation}\label{small-eig}
\lambda_1 = \frac{2 S a b P(+\infty)}{\int_{-\infty}^{\infty} K_{0y}^2 dy} \left( \tan^2(\sqrt{ab}) + 1 \right) = \frac{2 S a b P(+\infty)}{\int_{-\infty}^{\infty} K_{0y}^2 dy \cdot \cos^2(\sqrt{ab})}.
\end{equation}
Here, the constant $S$ is determined by the nonlinear matching condition \eqref{ssmatching}. Note that the right hand side of \eqref{P_Indentity} is negative, so we can deduce that $P(y)<0$. It follows $\lambda_1<0$.

\section{The slow drift dynamics of a single spike}\label{sec4}
When the spike is initially away from the origin, it starts to move slowly towards the origin, as illustrated in the Figure~ \ref{fig:intro(c)}. In this section, we analyze the slow dynamics of a one-spike solution to system \eqref{eq0}. We derive a differential equation that governs the motion of the spike's location, $x_0(t)$, defined as the point of maximum capital concentration. By linearizing this equation around the stable equilibrium $x_0 = 0$, we show that the decay rate of infinitesimal perturbations coincides with the small eigenvalue given in \eqref{small-eig}.

We assume the spike moves on a slow time scale, introducing
\begin{equation}
   T = \varepsilon^3 t, \quad x_0 = x_0(T).
\end{equation}
In the inner region near the spike, we define the stretched variable
\begin{equation}\label{inner_var}
    y = \frac{x - x_0(T)}{\varepsilon}, \quad L(y) = l(x_0 + \varepsilon y), \quad K(y) = k(x_0 + \varepsilon y),
\end{equation}
then \eqref{SS01} becomes
\begin{align}\label{SSinnerdynamics}
\begin{cases}
\tau \varepsilon^4 L_{y} \frac{\text{d} x_0 }{dT}   = L_{yy} - (L K_y)_y + \varepsilon^2 b L(a - L), & y \in \mathbb{R}, \\[4pt]
\varepsilon^2 K_{y} \frac{\text{d} x_0 }{dT} = K_{yy} - K + K^{\theta} L^{1-\theta}, & y \in \mathbb{R}, \\[4pt]
L_y(\pm \infty) = K_y(\pm \infty) = 0.
\end{cases}
\end{align}
We expand the inner solutions asymptotically as
\begin{equation}
L = L_0 + \varepsilon^2 L_1 + \cdots, \quad K = K_0 + \varepsilon^2 K_1 + \cdots.
\end{equation}
Substituting this expansion into \eqref{SSinnerdynamics} and collecting the leading order term, we obtain the system:
\begin{equation}\label{inner0_dyna}
\begin{aligned}
    0 &= L_{0yy} - (L_0 K_{0y})_y, & &y \in (-\infty, \infty), \\
    0 &= K_{0yy} - K_0 + K_0^{\theta} L_0^{1-\theta}, & &y \in (-\infty, \infty), \\
    &L_{0y}(\pm \infty) = K_{0y}(\pm \infty) = 0.
\end{aligned}
\end{equation}
We solve the first equation for $L_0$, imposing the far-field condition $L_0(\pm \infty) = S_0$, to find
\begin{equation}
    L_0 = S_0 e^{K_0 - S_0}.
\end{equation}
Here, the constant $S_0$ depends on the spike location $x_0$. Substituting this into the second equation yields an equation for $K_0$:
\begin{equation}\label{K0_dynamics}
    0 = K_{0yy} - K_0 + S_0^{1-\theta} K_0^{\theta} e^{(1-\theta)(K_{0} - S_0)}, \quad \text{with} \quad K_{0y}(0) = 0, \quad K_{0y}(\infty) = 0.
\end{equation}

At the next order $\mathcal{O}(\varepsilon^2)$, we obtain the system:
\begin{align}
    0 &= L_{1yy} - (L_1 K_{0y} + L_0 K_{1y})_y + b(a L_0 - L_0^2),  \label{L1dyna} \\
    -K_{0y} \frac{\dd x_0}{\dd T} &= K_{1yy} - K_1 + \theta K_0^{\theta-1} L_0^{1-\theta} K_1 + (1-\theta) K_0^{\theta} L_0^{-\theta} L_1,  \label{K1dyna}
\end{align}
with far-field conditions
\begin{equation}
L_{1y}(\pm \infty) = K_{1y}(\pm \infty), \quad L_{1}(\pm \infty) = K_{1}(\pm \infty).
\end{equation}
Unlike the steady-state analysis, we no longer require $K_1$ to be an even function. Integrating equation \eqref{L1dyna} over the real line yields
\begin{equation}
    \left[ L_{1y} + b(a S_0 - S_0^2) y \right]_{-\infty}^{\infty} = -b \int_{-\infty}^{\infty} \left[ (a L_0 - L_0^2) - (a S_0 - S_0^2) \right] dy.
\end{equation}

To derive the equation of motion for the spike, we apply a solvability condition. Specifically, we multiply \eqref{L1dyna} by $P$ and \eqref{K1dyna} by $Q$  (the components of the null vector of the adjoint operator solved in \eqref{QQ_S} and \eqref{PS}), add the resulting expressions, and integrate. This procedure yields the spike dynamics:
\begin{equation}\label{x_dynamics}
    \frac{\dd x_0}{\dd T}= -\frac{P(+\infty)}{\int_{-\infty}^{\infty} K_{0y}^2 dy} \left( L_{1y}(+\infty) + L_{1y}(-\infty) \right).
\end{equation}

\begin{figure}
    \centering
    \begin{minipage}{0.48\linewidth}
        \centering
        \includegraphics[width=\linewidth]{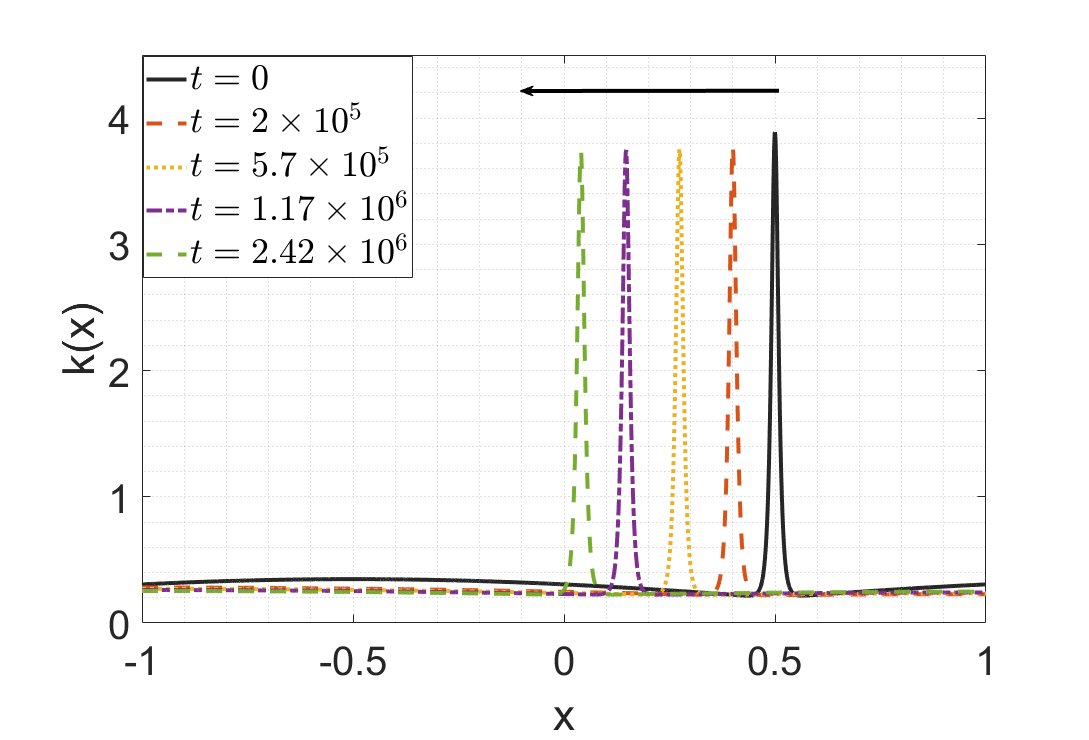}
        \subcaption{Snapshots of the spike profile at different times.}
    \end{minipage}
    \hfill
    \begin{minipage}{0.48\linewidth}
        \centering
        \includegraphics[width=\linewidth]{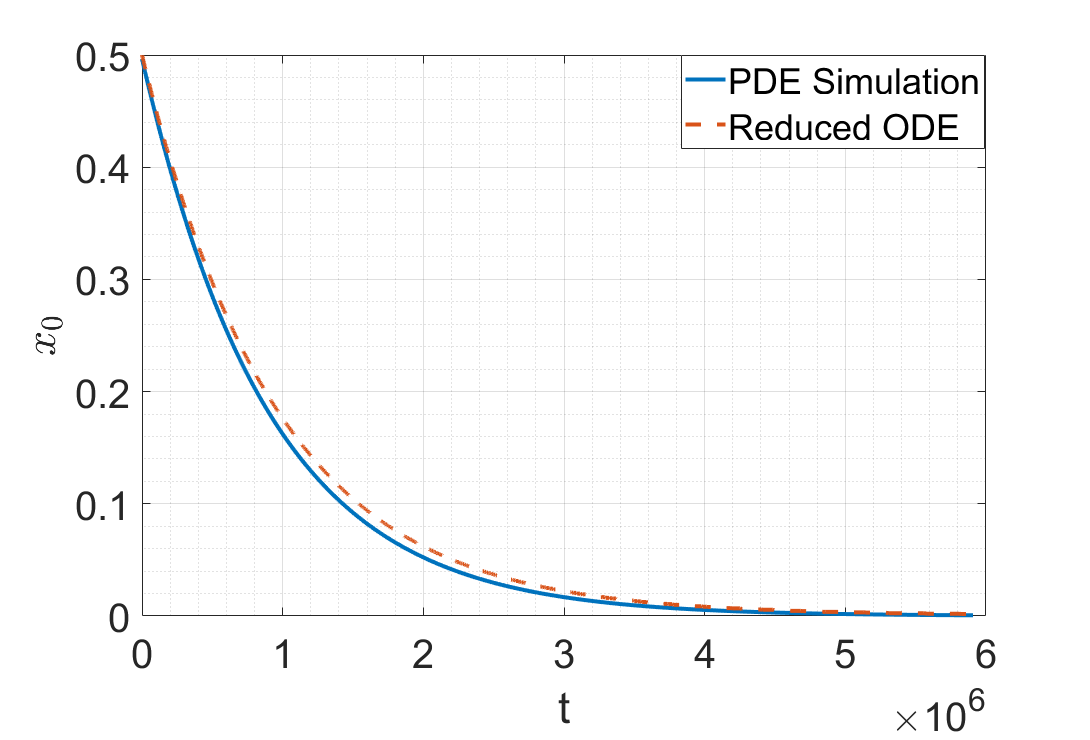}
        \subcaption{Trajectory of the spike center over time.}
    \end{minipage}
    \caption{Dynamics of a single interior spike in the PDE system \eqref{SS01}. 
    (a) Evolution of the spike profile, initialized at $x_0(0)=0.5$. 
    (b) Comparison of the spike center trajectory from direct PDE simulations (solid line) 
    with the reduced ODE prediction \eqref{x_dynamics_final} (dashed line). 
    Parameters: $\theta=0.5$, $\varepsilon=5\times10^{-3}$, $a=1$, $b=0.25$.}
    \label{fig:slowmotion}
\end{figure}

In the outer region, the equation for $l$ simplifies to
\begin{equation}
   l_{xx} - (l l_x)_x + b l(a-l) \sim 0  \quad \text{with} \quad l(x_0) = S_0, \quad l'(\pm 1) = 0.
\end{equation}
Since $S_0 \ll 1$, then
\begin{equation} \label{outer_appro_dynamics}
    l_{xx} + a b l \sim 0,\quad \text{with} \quad l(x_0) = S_0, \quad l'(\pm 1) = 0.
\end{equation}
The solution to this outer problem is
\begin{equation}\label{ldyna}
l(x) = 
    \begin{cases}
        S_0 \dfrac{\cos(\sqrt{ab}(x - 1))}{\cos(\sqrt{ab}(x_0 - 1))}, & x_0 < x < 1; \\[1em]
        S_0 \dfrac{\cos(\sqrt{ab}(x + 1))}{\cos(\sqrt{ab}(x_0 + 1))}, & -1 < x < x_0.
    \end{cases}
\end{equation}

Matching the derivative of this outer solution to the far-field behavior of the inner solution provides the condition
\begin{equation}\label{dynamics_match1}
\begin{aligned}
   -\varepsilon^2 b \int_{-\infty}^\infty \left[ a L_0 - L_0^2 - (a S_0 - S_0^2) \right] dy = \varepsilon \left( l_x(x_0^+) - l_x(x_0^-) \right)
\end{aligned}
\end{equation}
Simplifying \eqref{dynamics_match1} yields
\begin{equation}\label{matching_dynamics}
   2\varepsilon b\mathcal{I}(S_0)
   = - S_0 \sqrt{ab} \left( \tan(\sqrt{ab}(x_0 - 1)) - \tan(\sqrt{ab}(x_0 + 1)) \right).
\end{equation}
Equation \eqref{matching_dynamics} determines the dependence of the spike amplitude $S_0$ on its location $x_0$. Substituting the jump condition \eqref{matching_dynamics} in the outer derivative into \eqref{x_dynamics} gives the final equation of motion for the spike:
\begin{equation}\label{x_dynamics_final}
    \frac{\dd x_0}{\dd T} = \frac{ P(\infty)}{ \varepsilon \int_{-\infty}^{\infty} K_{0y}^2~dy } S_0 \sqrt{ab} \left( \tan(\sqrt{ab}(x_0 + 1)) + \tan(\sqrt{ab}(x_0 - 1)) \right),
\end{equation}
where
\begin{equation}
    P(\infty) = \frac{(1-\theta)}{2} \int_0^\infty \frac{1}{L_0(z)} \left( K_{0y}^2(z) - (K_0^2(z) - S_0^2) \right)~dz.
\end{equation}

We note that $x_0 = 0$ is an equilibrium point of \eqref{x_dynamics_final}. To study its stability, we linearize around this equilibrium by substituting $x_0 = \phi e^{\lambda_o t}$ and $S_0 = S$. A direct computation yields the linear growth rate:
\begin{equation}
    \lambda_o = \frac{\varepsilon^2 P(\infty)}{\int_{-\infty}^{\infty} K_{0y}^2~dy} \frac{2S a b}{\cos^2(\sqrt{ab})}.
\end{equation}
This result confirms that the dynamic decay rate of the spike displacement exactly coincides with the small eigenvalue $\lambda_1$ previously derived in \eqref{small-eig}.

\section{Discussion}

The spatial Solow model \eqref{tt1}, introduced by \cite{juchem2015capital}, exhibits remarkably rich dynamics.  Among the most striking behaviors is the emergence of aggregation patterns, or \emph{spikes}, in which capital and labor remain scarce across most of the domain but concentrate sharply in localized regions.  In this work, we focused on such aggregation 
phenomena under the asymptotic limit of small capital diffusion, with the capital-induced labor migration. Using the asymptotic–numerical hybrid method, we constructed a single-spike steady state and analyzed its stability via the linearized eigenvalue problem.  Our study reveals the existence of a stable single interior spike (a stable economic center).  In particular, when the reaction constant exceeds a critical threshold, the stability deteriorates, giving rise to oscillations in amplitude.  Furthermore, we derived a reduced ODE describing the slow drift of a spike as its center shifts away from the midpoint of the domain.  {Our mathematical analysis reveals that restricting capital mobility induces economic activity to concentrate within specific regions, resulting in the emergence of core–periphery structures characterized by interior spike profiles, as shown in Figure~\ref{fig:intro}.  This illustrates that higher capital mobility facilitates a more uniform distribution of economic growth across regions.}

\begin{figure}
    \centering
    \begin{subfigure}[b]{0.48\linewidth}
        \centering
        \includegraphics[width=\linewidth]{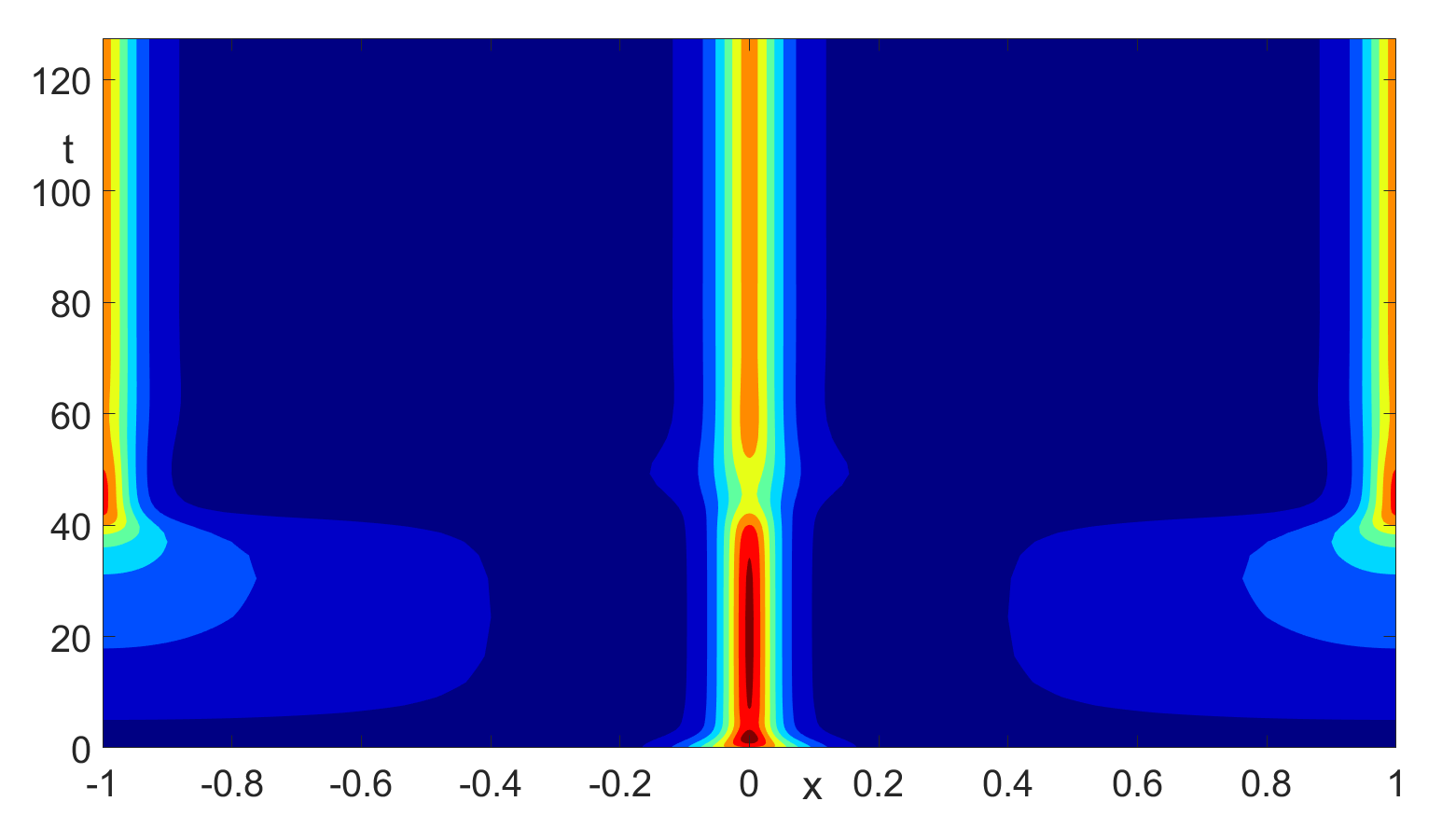}
        \caption{Spike insertion: two boundary spikes emerge after a single interior spike becomes unstable.}\label{fig:discussa}
    \end{subfigure}
    \hfill
    \begin{subfigure}[b]{0.48\linewidth}
        \centering
        \includegraphics[width=\linewidth]{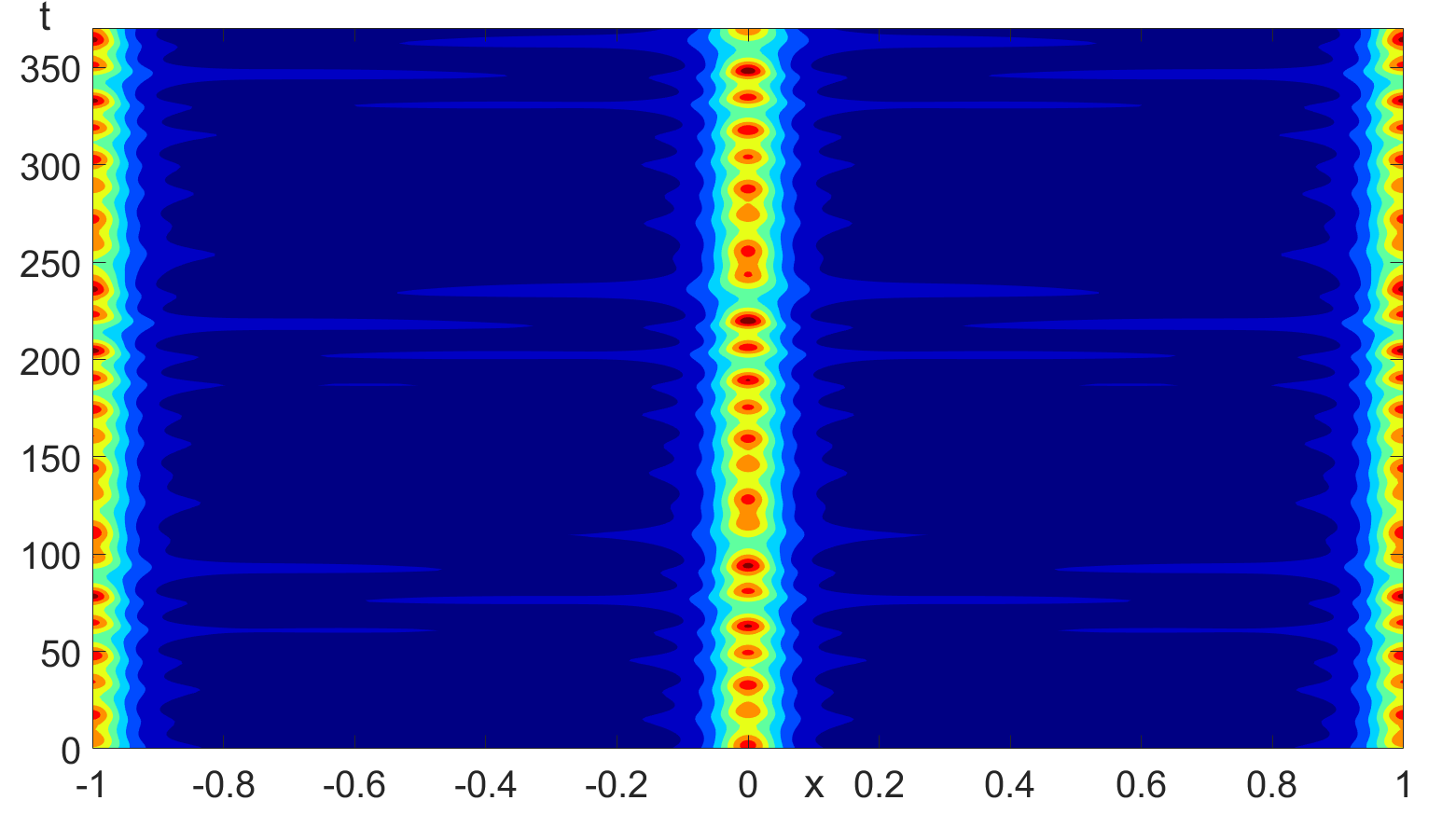}
        \caption{Blinking spikes: a interior spike and two boundary spikes oscillate in height in a complex manner.}
    \end{subfigure}
    \caption{Dynamics of a single spike and a single spike plus two boundary spikes in the spatial Solow model. Paramters are $a=2.5,~b=1,\varepsilon =2\times10^{-2}.$ (a) $~\tau=1$; (b) $~\tau=10$.  }
    \label{fig:discussion}
\end{figure}

For clarity, our analysis focuses solely on the case of a single interior spike. Nevertheless, the stability framework developed here can be naturally extended to multi-spike configurations via Floquet-type theory, as demonstrated in \cite{van2005stability,kong2024existence1}. Some interesting but yet unexplored profiles and dynamics are shown in Figure~\ref{fig:discussion}. We also note that the Keller-Segel model with logistic growth, studied in \cite{kong2024existence1}, corresponds to the special case of \eqref{SS01} with $\theta=0$.

There are several intriguing research directions that deserve the further exploration in the future.  First of all, it would be interesting to extend the model to incorporate the spatial distribution of technological effects \cite{capasso2010spatial,gonzalez2022mathematical,urena2024numerical}.
Incorporating the technology effect in (\ref{tt1}) would allow the model to capture more realistic economic scenarios such as technology-driven agglomeration, thereby providing a deeper understanding of the spatial dynamics of economic growth. 
 Also, the study of economic aggregation phenomena within the Solow model in higher spatial dimensions presents a promising direction for future research.  In higher dimensions, the geometry of the domain plays a crucial role in shaping the formation and stability of patterns \cite{tzou2017stability,xie2017moving,tzou2023oscillatory,deng2024oscillatory}.  Similarly to patterns observed in the Keller-Segel model, higher-dimensional settings can give rise to more intricate spatial structures, including spots \cite{kong2023existence2}, rings \cite{chen2020stationary}, stripes \cite{jin2016pattern}, and mesa-like configurations \cite{carrillo2020phase}, which are not observed in one-dimensional settings. Studying these higher-dimensional patterns could provide valuable insights into how spatial constraints, boundary conditions, and multi-dimensional interactions shape the distribution of capital and labor across regions, potentially revealing new mechanisms of economic clustering and agglomeration.

\appendix
\begin{section}{Subinner approximation of the NLEP problem} \label{app:A}
To investigate NLEP \eqref{Nlep_1} further, we analyze \eqref{Nlep_1} using a sub-inner region approximation. In the sub-inner region, where $y\sim    \mathcal{O}(\xi^{-1})$ and $K_0=\xi\gg 1$, we introduce the rescaled variable
\[
z=\frac{\xi(1-\theta)y}{2}.
\]
At the order of $\mathcal{O}(\xi^2)$, using the approximation  $\varepsilon b 
    \sim 2S\frac{\sqrt{ab}\tan(\sqrt{ab}}{\int_{-\infty}^\infty L_0^2~dy }  $ from \eqref{subalge},  equation \eqref{Nlep_1} reduces to
\begin{align}\label{psi00}
    \lambda_0\Psi_{00} &\sim \frac{(1-\theta)^2\xi^2}{4} \Psi_{00zz} \\
    &+ (1-\theta) \xi^\theta S^{1-\theta} e^{(1-\theta)\xi}e^{K_{00}(z)}\left( \Psi_{00}- \frac{2\sqrt{ab}\tan(\sqrt{ab})}{2\sqrt{ab}\tan(\sqrt {ab})-\mu\tan\mu} \frac{\int_{-\infty}^{+\infty} L_{00}^2\Psi_{00}\,dz}{\int_{-\infty}^{+\infty}L_{00}^2\,dz} \right).\nonumber
\end{align}
Using \eqref{subalge}, we obtain
\begin{align*}
   (\lambda_0+1)\Psi_{00} &\sim \frac{(1-\theta)^2\xi^2}{4} \Psi_{00zz} \\
   &+ \frac{\xi^2(1-\theta)^2}{2}e^{K_{00}(z)} \left( \Psi_{00}-\frac{2\sqrt{ab}\tan(\sqrt ab)}{2\sqrt{ab}\tan(\sqrt {ab})-\mu\tan\mu}\frac{\int_{-\infty}^{+\infty} L_{00}^2\Psi_{00}\,dz}{\int_{-\infty}^{+\infty}L_{00}^2\,dz} \right),
\end{align*}
or equivalently
\begin{equation}\label{psi000}
    \frac{4(\lambda_0+1)}{(1-\theta)\xi^2} \Psi_{00} \sim   \Psi_{00zz} + 2 e^{K_{00}(z)} \left( \Psi_{00}-\frac{2\sqrt{ab}\tan(\sqrt ab)}{2\sqrt{ab}\tan(\sqrt {ab})-\mu\tan\mu}\frac{\int_{-\infty}^{+\infty} L_{00}^2\Psi_{00}\,dz}{\int_{-\infty}^{+\infty}L_{00}^2\,dz} \right).
\end{equation}
Letting $U=2\sech^2(z)$, equation \eqref{psi000} becomes
\begin{equation}\label{psiNLEP}
   \frac{4(\lambda_0+1)}{(1-\theta)\xi^2} \Psi_{00} \sim   \Psi_{00zz} + U \left(\Psi_{00}-\frac{2\sqrt{ab}\tan(\sqrt ab)}{2\sqrt{ab}\tan(\sqrt {ab})-\mu\tan\mu}\frac{\int_{-\infty}^{+\infty} U^{\frac{2}{1-\theta}}\Psi_{00}\,dz}{\int_{-\infty}^{+\infty}U^{\frac{2}{1-\theta}}\,dz} \right).
\end{equation}

This motivates the study of the nonlocal eigenvalue problem
\begin{equation}\label{nlepv0}
    \Psi_{zz} + U\Psi -\alpha U \frac{\int_{-\infty}^{+\infty} U^{2r} \Psi\,dz}{\int_{-\infty}^{+\infty} U^{2r}\,dz}= \Lambda \Psi, \qquad r\geq1.
\end{equation}
then $\lambda_0$ and $\Lambda$ have the same sign when $\xi$ is large, since
\begin{equation}
    \lambda_0=\frac{1}{4}\Lambda \xi^2(1-\theta)-1 \sim \frac{1}{4}\Lambda \xi^2(1-\theta)\qquad \xi \to \infty. 
\end{equation}
Upon introducing $\bar z=2z$, \eqref{nlepv0} becomes
\begin{equation}\label{nlepv1}
    \Psi_{\bar z \bar z} + \frac{\sech^2(\tfrac{\bar z}{2})}{2} \Psi -\alpha \frac{\sech^2(\tfrac{\bar z}{2})}{2} \frac{\int_{-\infty}^{+\infty} U^{2r} \Psi\,dz}{\int_{-\infty}^{+\infty} U^{2r}\,dz}= \Lambda \Psi, \qquad r\geq 1.
\end{equation}
Defining
\[
\rho=\frac{3}{2} \sech^2\!\left(\tfrac{\bar z}{2}\right),
\]
we can rewrite \eqref{nlepv1} as
\begin{equation}\label{nlepv2}
    \Psi_{\bar z \bar z} + \tfrac{\rho}{3} \Psi - \tfrac{\alpha}{3}\frac{\int_{-\infty}^{+\infty} \rho^{2r} \Psi\,dz}{\int_{-\infty}^{+\infty} \rho^{2r}\,dz}\rho= \Lambda \Psi, \qquad r\geq 1.
\end{equation}
Since $\rho$ satisfies the nonlinear ODE
\[
\rho''-\rho+\rho^2=0,
\]
the problem \eqref{nlepv2} falls into the class of NLEPs analyzed in \cite{kolokolnikov2009existence} (see Theorem~3 therein), where they have derived the following lemma:
\begin{lemma}\label{lemma}
The nonlocal eigenvalue problem \eqref{nlepv2} admits a positive eigenvalue if $\alpha<1$. For the case $r=1$, if $\alpha>1$, the spectrum contains no eigenvalue with positive real part.
\end{lemma}

It is worth noting that  the authors 
conjectured in \cite{kolokolnikov2009existence} that Lemma \eqref{lemma} remains valid even for $r>1.$  In our case, $\alpha=\frac{2\sqrt{ab}\tan(\sqrt ab)}{2\sqrt{ab}\tan(\sqrt {ab})-\mu\tan\mu} =2 >1$ when $\tau=0$.  If the conjecture holds, then the problem \eqref{psiNLEP} admits no eigenvalues with positive real parts, implying that the single-spike profile is stable when $\tau=0$ in the limit $\xi\gg 1.$ 
\end{section}

\bibliographystyle{unsrt}
\bibliography{main}

\end{document}